\documentclass[final,1p,times]{elsarticle}

\usepackage{amssymb}
 \usepackage{amsthm}
\usepackage{amscd}
\usepackage{bm}
\usepackage{amsmath}
\usepackage{amsfonts}
\usepackage{amssymb}
\usepackage{graphicx}
\newtheorem{theorem}{Theorem}

\newtheorem{lemma}[theorem]{Lemma}

\usepackage{mathrsfs}
\usepackage{titletoc}

\renewcommand{\d}{\delta}

\newcommand{\e}{\epsilon}
\newcommand{\ra}{\rightarrow}
\newcommand{\p}{\partial}
\newcommand{\f}{\frac}

\newcommand{\be}{\begin{equation}}
\renewcommand{\ra}{\rightarrow}
\newcommand{\ee}{\end{equation}}
\newcommand{\bea}{\begin{eqnarray}}
\newcommand{\eea}{\end{eqnarray}}
\newcommand{\bna}{\begin{eqnarray*}}
\newcommand{\ena}{\end{eqnarray*}}

\renewcommand{\le}{\left}
\newcommand{\ri}{\right}

\newcommand{\Si}{\Sigma}

\journal{***}

\begin{document}

\begin{frontmatter}

\title{Trudinger-Moser inequalities on a closed Riemann surface with a symmetric conical metric}

\author{Yu Fang}
\ead{fangyu-3066@ruc.edu.cn }

\author{Yunyan Yang\footnote{corresponding author}}
 \ead{yunyanyang@ruc.edu.cn}
 
\address{ Department of Mathematics,
Renmin University of China, Beijing 100872, P. R. China}

\begin{abstract}
This is a continuation of our previous work \cite{preprint}.
Let $(\Sigma,g)$ be a closed Riemann surface, where the metric $g$ has conical singularities at finite points.
Suppose $\mathbf{G}$ is a group whose elements are isometries acting on $(\Sigma,g)$. Trudinger-Moser inequalities
involving $\mathbf{G}$ are established via the method of blow-up analysis, and the corresponding extremals are also
obtained. This extends previous results of Chen \cite{W. Chen}, Iula-Manicini \cite{macini},
and the authors \cite{preprint}.
  \end{abstract}


\end{frontmatter}

\titlecontents{section}[0mm]
                       {\vspace{.2\baselineskip}}%\bfseries}
                       {\thecontentslabel~\hspace{.5em}}
                        {}
                        {\dotfill\contentspage[{\makebox[0pt][r]{\thecontentspage}}]}
\titlecontents{subsection}[3mm]
                       {\vspace{.2\baselineskip}}%\bfseries}
                       {\thecontentslabel~\hspace{.5em}}
                        {}
                       {\dotfill\contentspage[{\makebox[0pt][r]{\thecontentspage}}]}

\setcounter{tocdepth}{2}
%\tableofcontents

%\setcounter{tocdepth}{1}
%\tableofcontents

%\tableofcontents
\section{Introduction and main results}\label{introduction}
Let $\mathbb{S}^2$ be the $2$-dimensional sphere $x_1^2+x_2^2+x_3^2=1$ endowed with a metric $g_1=dx_1^2+dx_2^2+dx_3^2$ for all $x=(x_1,x_2,x_3)\in\mathbb{R}^3$. It was proved by Moser \cite{Moser} that there exists a universal constant $C$ satisfying
\begin{equation}\label{M1}
\int_{\mathbb{S}^2}e^{4\pi u^{2}}dv_{g_1}\leq C
\end{equation}
for all smooth functions $u$ with $\int_{\mathbb{S}^2}|\nabla_{g_1}u|^2 dv_{g_1}\leq 1$ and $\int_{\mathbb{S}^2}u dv_{g_1}=0$,
where $\nabla_{g_1}$ and $dv_{g_1}$ stand for the gradient operator and the volume element on $(\mathbb{S}^2,g_1)$ respectively.
Here $4\pi$ is best constant in the sense that when $4\pi$ is replaced by any $\alpha>4\pi$,
the integrals are still finite, but the universal constant $C$ no longer exists.
It was also remarked by Moser \cite{Jmoser} that if one considers even functions $u$, say $u(x)=u(-x)$ for all $x\in\mathbb{S}^2$,
then the constant $4\pi$ in (\ref{M1}) would double. Namely there exists an absolute constant $C$ such that
  \begin{equation}\label{1.3}
\int_{\mathbb{S}^2}e^{8\pi u^{2}}dv_{g_1}\leq C
\end{equation}
for all even functions $u$ satisfying $\int_{\mathbb{S}^2}|\nabla_{g_1}u|^2 dv_{g_1}\leq 1$,
$\int_{\mathbb{S}^2}u dv_{g_1}=0$.

A general manifold version of (\ref{M1}) was established by Fontana \cite{Fontana} via the estimation on Green functions and O'Neil's lemma \cite{ONeil}.
This comes from an Euclidean scheme designed by  Adams \cite{Adams}. However,  Li
\cite{Lijpde} was able to prove the inequality (\ref{M1}) by the method of blow-up analysis. In a recent work \cite{preprint}, we extended (\ref{1.3}) to the case of closed Riemann surface
with a smooth ``symmetric" metric.
In the current paper, we consider the case of closed Riemann surface with a ``symmetric" singular  metric. For earlier works
on Trudinger-Moser inequalities involving singular metrics, we refer the reader to Troyanov \cite{Troyanov}, Chen \cite{W. Chen}, Adimurthi-Sandeep \cite{A-S}, Adimurthi-Yang \cite{Adi-Yang}, Li-Yang \cite{Li-Yang-JDE}, Csato-Roy \cite{Csato}, Yang-Zhu \cite{2017JFA}, Iula-Mancini \cite{macini}
and the references therein.

Now we recall some notations from differential geometry. Let $(\small{\Si},g_0)$ be a closed Riemann surface, and $d_{g_0}(\cdot,\cdot)$ be the geodesic distant between two points of $\Si$. A  smooth metric $g$ defined on $\Si\setminus\{p_1,\cdots,p_{L}\}$ is said to have conical singularity of order $\beta_i>-1$ at $p_i$, $i=1,\cdots,L$, if
\be\label{rho}
g=\rho g_0,
\ee where $\rho\in C^\infty(\Si\setminus\{p_1,...,p_{L}\},g_0)$ satisfies $\rho>0$ on $\Si\setminus\{p_1,...,p_{L}\}$ and
\be\label{h1}
0<C\leq\f{\rho(x)}{d_{g_0}(x,p_i)^{2\beta_i}}\in C^0(\Sigma,g_0)\ee
for some constant $C$ and  $i=1,\cdots,L$. Here we write the righthand side of (\ref{h1}) in the sense that $\rho/d_{g_0}(x,p_i)^{2\beta_i}$ can be continuously extended to
the whole surface $(\Sigma,g_0)$.
 With (\ref{rho}) and (\ref{h1}),  $(\Sigma,g)$ is called a closed Riemann surface having conical singularities of the divisor $\mathbf{b}=\sum_{i=1}^L\beta_ip_i$.
For more details on singular surface, we refer the reader to Troyanov \cite{Troyanov}.
We say that $\mathbf{G}=\{\sigma_1,\sigma_2,\cdots,\sigma_{N}\}$ is a finite isometric group acting on ($\Sigma,g)$, if
 each smooth map $\sigma_k: \Si\ra\Si$
satisfies
\be\label{isorho}
(\sigma_k^*{g_0})_x={g_0}_{\sigma_k(x)} \quad{\rm and}\quad\rho(\sigma_k(x))=\rho(x)\quad{\rm for\,\,\,all}\quad x\in\Sigma.
\ee
This in particular implies
\be\label{iso}
\sigma^*g_x=g_{\sigma(x)}\quad{\rm for\,\,\,all}\quad x\in\Sigma.
\ee
Note that $\mathbf{G}$ is a geometric structure on special Riemann surface $(\Sigma,g)$.
It is clear that $\mathbf{G}(p_j)=\{\sigma_i(p_j)\}_{i=1}^N\subset\{p_1,\cdots,p_L\}$
for all $j$,
and that $\beta_k=\beta_j$ provided that $p_k\in\mathbf{G}(p_j)$ for some $j$. Denote for any $x\in\Sigma$,
\be\label{number}I(x)=\sharp \mathbf{G}(x)\ee
and
\be\label{betax}
\beta(x)=\le\{
  \begin{array}{lll}
0,\quad x\not\in \{p_1,\cdots,p_L\},\\[1.5ex]
\beta_j, \quad x=p_j,\,\,1\leq j\leq L,
  \end{array}
  \ri.
\ee
where $\sharp \mathbf{A}$ is the number of all distinct elements in the set $\mathbf{A}$. Noting that $1\leq I(x)\leq N$
and $\beta(x)>-1$ for all $x\in\Sigma$, one defines
\be\label{M0}
\ell=\min_{x\in\Si}\min\le\{I(x),I(x)(1+\beta(x))\ri\}.
\ee

Let $W^{1,2}(\Si,g)$ be the completion of $C^\infty(\Sigma,g_0)$  under the norm
\be\label{w12}
\|u\|_{W^{1,2}(\Si,g)}=\le(\int_\Si \le(|\nabla_g u|^2+u^2 \ri)dv_g \ri)^{1/2}.
\ee
 For convenience, a subspace of $W^{1,2}(\Si,g)$ is denoted by
\bea
\label{hspace}
\mathscr{H}_{\mathbf{G}}=\left\{u\in W^{1,2}(\Sigma,g):\int_\Si u dv_g=0,\,\,u(x)=u(\sigma(x))\,\,{\text{for}}\,\,\text{a.e}.\,\, x\in
\Si  \,\,\text{and\, all}\,\,\sigma\in\mathbf{G}\right\}.
\eea
Clearly, $\mathscr{H}_{\mathbf{G}}$ is a Hilbert space with an inner product
 $$\langle u,v\rangle_{\mathscr{H}_{\mathbf{G}}}=\int_\Si \langle\nabla_g u,\nabla_g v\rangle dv_g.$$
The first eigenvalue of $\Delta_g$ on $\mathscr{H}_{\mathbf{G}}$ reads
\be\label{lamp}
\lambda_1^{\mathbf{G}}=\inf_{u\in\mathscr{H}_{\mathbf{G}},\,\int_\Si u^2 dv_g=1}\int_{\Si}|\nabla_{g}u|^{2}\mathrm{d}v_{g},
\ee
where $\Delta_g$ is the Laplace-Beltrami operator with respect to the conical metric $g$.
A direct method of variation leads to $\lambda_1^{\mathbf{G}}>0$. For any $\alpha$ strictly less than $\lambda_1^{\mathbf{G}}$,
we can define an equivalent norm of (\ref{w12}) on $\mathscr{H}_{\mathbf{G}}$ by
\be\label{a,2}
\|u\|_{1,\alpha}=\le(\int_\Si |\nabla_g u|^2dv_g-\alpha\int_\Si u^2dv_g\ri)^{1/2}.
\ee
The first eigenfunction space with respect to $\lambda_1^{\mathbf{G}}$ reads as
\be\label{1-stspace}
E_{\lambda_1^{\mathbf{G}}}=\le\{u\in\mathscr{H}_{\mathbf{G}}: \Delta_g u=\lambda_1^{\mathbf{G}}u\ri\}.
\ee
According to Chen \cite{W. Chen}, there holds
\bea\label{chenine}
\sup\limits_{u\in \mathscr{H}_{\mathbf{G}},\int_\Sigma|\nabla_g u|^2dv_g\leq1}  \int_{\Si}e^{4\pi \ell u^2} dv_g<\infty,
\eea
where $\ell$ is given as in (\ref{M0}), and $4\pi\ell$ is the best constant for (\ref{chenine}) in the sense that if $4\pi\ell$
is replaced by any $\gamma>4\pi\ell$, then the supremum in (\ref{chenine}) is infinity.
We first concern the attainability of the above supremum and have the following more general result:
\begin{theorem}\label{TH1} Let $(\Sigma,g)$ be a closed Riemann surface with conical singularities of the divisor $\mathbf{b}=\sum_{i=1}^L\beta_ip_i$, where $p_i$ belongs to $\Sigma$ and
\be\label{assumption}{-1}<\beta_i\leq 0,\quad i=1,\cdots,L.\ee
Suppose that $\mathbf{G}=\{\sigma_1,\sigma_2,\cdots,\sigma_{N}\}$ is a group of isometries given as in (\ref{isorho}), and that $\ell$, $\mathscr{H}_{\mathbf{G}}$ and $\lambda_{1}^{\mathbf{G}}$ are defined as in (\ref{M0}), (\ref{hspace}) and (\ref{lamp})
respectively. Then for any  $\alpha<\lambda_1^{\mathbf{G}}$, the supremum
\be\label{Th1}
\sup\limits_{ u\in \mathscr{H}_{\mathbf{G}},\,\| u\|_{1,\alpha}\leq 1}\int_{\Sigma}e^{4\pi \ell u^2} dv_g
\ee
is attained by some function $u_{0}\in C^{1}(\Sigma\setminus\{p_1,\cdots,p_L\},g_0)\cap C^0(\Sigma,g_0)\cap\mathscr{H}_{\mathbf{G}}$
satisfying $\|u_{0}\|_{1,\alpha}=1$, where $g_0$ is  a smooth metric given as in (\ref{rho}) and $\|\cdot\|_{1,\alpha}$ is defined as in (\ref{a,2}).
\end{theorem}
When $N=1$, Theorem \ref{TH1} reduces to one of results of Iula-Mancini \cite{macini}. While if $\beta(x)\equiv 0$ for all $x\in \Sigma$,
 then Theorem \ref{TH1} is exactly our earlier result \cite{preprint}.
  To prove Theorem \ref{TH1}, we use the method of blow-up analysis designed by Li \cite{Lijpde}. Early groundbreaking
works go back to Carleson-Chang \cite{CC}, Ding-Jost-Li-Wang \cite{DJLW} and Adimurthi-Struwe \cite{Adi-Stru}.

As in our previous work \cite[Theorem\,2]{preprint}, we may also consider the effect of higher order eigenvalues of $\Delta_g$ on  Trudinger-Moser inequalities. Set $E_0=\{0\}$, $E_0^{\bot}=\mathscr{H}_{\mathbf{G}}$, and $E_1=E_{\lambda_1^{\mathbf{G}}}$ is defined as in (\ref{1-stspace}).
By induction, $E_j$ and $E_j^\perp$ can be defined for any positive integer $j$.
 In precise, for any $j\geq 1$, we set
$E_{j}=E_{\lambda_1^{\mathbf{G}}}\oplus \cdots \oplus E_{\lambda_{j}^{\mathbf{G}}}$ and
\begin{equation}\label{elbot}
 E_{j}^{\bot}=\le\{u\in \mathscr{H}_{\mathbf{G}}:\int_{\Si} uv dv_g=0,\;\forall v\in E_{j} \ri\},
\end{equation}
where $\lambda_{j}^{\mathbf{G}}$ is the $j$-th eigenvalue of $\Delta_g$ written by
\be\label{laml}
\lambda_{j}^{\mathbf{G}}=\inf_{u\in E_{j-1}^{\bot},\,\int_\Si u^2 dv_g=1}\int_{\Si}|\nabla_{g}u|^{2}\mathrm{d}v_{g},
\ee
and $E_{\lambda_{j}^{\mathbf{G}}}=\{u\in E_{j-1}^{\bot}: \Delta_g u=\lambda_{j}^{\mathbf{G}}u\}$ is the corresponding $j$-th eigenfunction space.
Obviously for any fixed $\alpha<\lambda_{j+1}^{\mathbf{G}}$, $\|\cdot\|_{1,\alpha}$ is equivalent to $\|\cdot\|_{W^{1,2}(\Si,g)}$ on the space $E_{j}^{\bot}$.\\

Our second result reads as follows:
\begin{theorem}\label{TH2}
Let $(\Sigma,g)$ be a closed Riemann surface with conical singularities of divisor $\mathbf{b}=\sum_{i=1}^L\beta_ip_i$, where $p_i$ belongs to $\Sigma$ and ${-1}<\beta_i\leq 0$ for $i=1,\cdots,L$. Suppose that $\mathbf{G}=\{\sigma_1,\sigma_2,\cdots,\sigma_{N}\}$ is a group of isometries given as in (\ref{isorho}).
Then for any integer $j\geq 1$ and any real number $\alpha$ satisfying $\alpha<\lambda_{j+1}^{\mathbf{G}}$, the supremum
\be\label{Th2}
\sup\limits_{ u\in E_{j}^{\bot},\,\| u\|_{1,\alpha}\leq 1}\int_{\Sigma}e^{4\pi \ell u^2} dv_g
\ee
can be attained by some function $u_{0}\in C^{1}(\Sigma\setminus\{p_1,\cdots,p_L\},g_0)\cap C^0(\Sigma,g_0)\cap E_{j}^{\bot}$ with $\|u_{0}\|_{1,\alpha}=1$, where
$\lambda_{j+1}^{\mathbf{G}}$, $E_{j}^{\bot}$, $\ell$ and $\|\cdot\|_{1,\alpha}$ are defined as in (\ref{lamp}), (\ref{elbot}), (\ref{M0}), and (\ref{a,2}) respectively, and $g_0$ is a smooth metric given as in (\ref{rho}).
\end{theorem}
The proof of Theorem \ref{TH2} is similar to that of Theorem \ref{TH1}. The difference is that we work on the space $E_j^\perp$
instead of $\mathscr{H}_{\mathbf{G}}$. Note that $E_j^\perp$ is still a Hilbert space for any $j\geq 1$.
For more details of Trudinger-Moser inequalities involving eigenvalues, we refer
 the reader to \cite{Adimurthi-D, Tintarev, Yang-JDE2015}.
In both proofs of Theorems \ref{TH1} and \ref{TH2}, to derive an upper bound of the Trudinger-Moser functional, we need a singular version of Carleson-Chang's estimate,
which was in literature due to Csato-Roy \cite{Csato} (see also Iula-Mancini \cite{macini}
and Li-Yang \cite{Li-Yang-JDE}), namely
\begin{lemma}\label{circle}\,Let $\mathbb{B}_r\subset\mathbb{R}^2$ be a ball centered at $0$ with radius $r$.
 If $\phi_\e\in W_0^{1,2}(\mathbb{B}_r)$ satisfies $\int_{\mathbb{B}_r}|\nabla \phi_\epsilon|^2dx\leq 1$, and $\phi_\e\rightharpoonup 0\,\text{weakly \,in}\,W_0^{1,2}(\mathbb{B}_r)$, then for any $\beta$ with $-1<\beta\leq0$, there holds
\be\label{cclemma}
\limsup_{\e\ra0}\int_{\mathbb{B}_r}e^{(1+\beta)4\pi  \phi_\e^2}|x|^{2\beta} dx\leq\int_{{\mathbb{B}_r}}|x|^{2\beta}dx+\frac{\pi e }{1+\beta}r^{2+2\beta}.
\ee
\end{lemma}
The proof of Lemma \ref{circle} is based on a rearrangement argument, Hardy-Littlewood inequality, and Carleson-Chang's estimate \cite{CC}. In the remaining part of this paper, we prove Theorems \ref{TH1} and \ref{TH2} in Sections \ref{proof} and \ref{pf-2} respectively. Throughout
this paper, we do not distinguish sequence and subsequence. Constants are often denoted by the same $C$ from line to line, even on the same line.

\section{Trudinger-Moser inequalities involving the first eigenvalue}\label{proof}
In this section we shall prove Theorem \ref{TH1} by using the method of blow-up analysis, which was originally used in this topic by Li \cite{Lijpde,LiScience},  and extensively used by
Yang \cite{Yang-Tran2007,Yang-JDE2015}, Li-Yang \cite{Li-Yang-JDE}, de Souza-do O \cite{de Souza-do O}, Yang-Zhu \cite{2017JFA},
Iula-Mancini \cite{macini} and others. The proof is divided into several subsections below.

\subsection{The best constant}
 Let $\ell$ be defined as in (\ref{M0}). It was proved by Chen \cite{W. Chen} that
\be\label{Chen}\sup_{u\in \mathscr{H}_{\mathbf{G}},\,\int_{\Si}|\nabla u|^2dv_g\leq1}\int_{\Si}{e^{\gamma u^2}} dv_g<\infty,\quad\forall
\gamma\leq 4\pi\ell;\ee
moreover,  the above integrals are still finite for any $\gamma>4\pi\ell$, but the supremum
\be\label{best}\sup_{u\in \mathscr{H}_{\mathbf{G}},\,\int_{\Si}|\nabla u|^2dv_g\leq1}\int_{\Si}{e^{\gamma u^2}} dv_g=\infty, \quad\forall \gamma>4\pi\ell.\ee
We now take the first eigenvalue $\lambda_1^{\bm G}$ of $\Delta_g$ (see (\ref{lamp}) above)  into account and have the following:
\begin{lemma}\label{initial}
For any $\alpha<\lambda_1^{\mathbf{G}}$, there exists a real number $\gamma_0>0$ such that
$$\sup_{u\in \mathscr{H}_{\mathbf{G}},\,\|u\|_{1,\alpha}\leq1}\int_{\Si}{e^{\gamma_0 u^2}} dv_g<\infty,$$
where $\|\cdot\|_{1,\alpha}$ is defined as in (\ref{a,2}).
\end{lemma}
\proof Assume $\alpha<\lambda_1^{\bm G}$ and $\|u\|_{1,\alpha}\leq 1$. Then
$$\le(1-\f{\alpha}{\lambda_1^{\mathbf{G}}}\ri)\int_\Sigma|\nabla_gu|^2dv_g\leq \int_\Sigma|\nabla_gu|^2dv_g-\alpha\int_\Sigma u^2dv_g\leq1.$$
This together with (\ref{Chen}) implies the existence of $\gamma_0$, as desired. $\hfill\Box$ \\

In view of Lemma \ref{initial}, for any fixed $\alpha<\lambda_1^{\mathbf{G}}$, we set
$$\gamma^\ast=\sup\le\{\gamma_0: \sup_{u\in \mathscr{H}_{\mathbf{G}},\,\|u\|_{1,\alpha}\leq1}\int_{\Si}{e^{\gamma_0 u^2}} dv_g<\infty\ri\}.$$
\begin{lemma}\label{bestconstant}
There holds $\gamma^\ast\geq 4\pi\ell$.
\end{lemma}

\proof Suppose $\gamma^\ast<4\pi\ell$. Then there exists a real number $\gamma_1$ with $\gamma^\ast<\gamma_1<4\pi\ell$
and a function sequence $(u_j)\subset \mathscr{H}_{\mathbf{G}}$
such that $\|u_j\|_{1,\alpha}\leq 1$ and
\be\label{infin}\int_\Sigma e^{\gamma_1 u_j^2}dv_g\ra\infty\quad{\rm as}\quad j\ra\infty.\ee
Since $\alpha<\lambda_1^{\mathbf{G}}$, we have
that $(u_j)$ is bounded in $W^{1,2}(\Sigma,g)$. Thus, $u_j$ converges to some $u_0$ weakly in
$W^{1,2}(\Sigma,g)$, strongly in $L^2(\Sigma,g)$ and almost everywhere in $\Sigma$. This particularly leads to
$$\|u_j-u_0\|_{1,\alpha}^2=\|u_j\|_{1,\alpha}^2-\|u_0\|_{1,\alpha}^2+o_j(1).$$
Clearly $u_0\in\mathscr{H}_{\mathbf{G}}$. We now claim that $u_0\equiv 0$. For otherwise, since $\|u_j\|_{1,\alpha}\leq 1$, there must hold
\be\label{u0}\int_\Sigma|\nabla_g(u_j-u_0)|^2dv_g\leq 1-\frac{1}{2}\|u_0\|_{1,\alpha}^2\ee
for sufficiently large $j$. Noting that $u_j^2\leq (1+\nu)(u_j-u_0)^2+(1+\nu^{-1})u_0^2$ for any $\nu>0$,
and that $e^{u_0^2}\in L^q(\Sigma,g)$ for all $q>1$, we conclude from (\ref{Chen}) and (\ref{u0}),
\be\label{leq}\int_\Sigma e^{\gamma_1 u_j^2}dv_g\leq C\ee
for some constant $C$ depending only on $\gamma_1$, $\ell$ and $u_0$. This contradicts (\ref{infin}) and confirms our claim
$u_0\equiv 0$. As a consequence
$$\int_\Sigma|\nabla_gu_j|^2dv_g\leq 1+\alpha\int_\Sigma u_j^2dv_g=1+o_j(1).$$
This together with (\ref{Chen}) gives (\ref{leq}), which again contradicts (\ref{infin}) and thus completes the proof of the lemma.
$\hfill\Box$\\

More precisely we have

\begin{lemma}\label{b-2} There holds
$\gamma^\ast=4\pi\ell$.
\end{lemma}
\proof By Lemma \ref{bestconstant}, $\gamma^\ast\geq 4\pi\ell$. Suppose $\gamma^\ast>4\pi\ell$. Fix some $\gamma_2$ with
$4\pi\ell<\gamma_2<\gamma^\ast$.  In view of (\ref{best}), there exists a function sequence $(M_k)\subset\mathscr{H}_{\mathbf{G}}$
such that \be\label{eng-1}\int_\Sigma|\nabla_gM_k|^2dv_g\leq 1\ee
and
\be\label{eng-2}\int_\Sigma e^{\gamma_2 M_k^2}dv_g\ra \infty.\ee
Obviously $(M_k)$ is bounded in $W^{1,2}(\Sigma,g)$. With no loss of generality, we assume
$M_k$ converges to $M_0$ weakly in $W^{1,2}(\Sigma,g)$, strongly in $L^2(\Sigma,g)$, and almost everywhere in $\Sigma$. Using the same
argument as in the proof of Lemma \ref{bestconstant}, we have $M_0\equiv 0$. It then follows that
\be\label{b-3}\|M_k\|_{1,\alpha}^2=\int_\Sigma|\nabla_gM_k|^2dv_g-\alpha\int_\Sigma M_k^2dv_g=1+o_k(1).\ee
Combining (\ref{eng-2}) and (\ref{b-3}), we have for some $\gamma_3$ with $\gamma_2<\gamma_3<\gamma^\ast$,
$$\sup_{u\in \mathscr{H}_{\mathbf{G}},\,\|u\|_{1,\alpha}\leq1}\int_{\Si}{e^{\gamma_3 u^2}} dv_g=\infty.$$
This contradicts the definition of $\gamma^\ast$. Therefore $\gamma^\ast$ must be $4\pi\ell$. $\hfill\Box$

\subsection{Maximizers for subcritical functionals}

In this subsection, using a direct method of variation, we show existence of maximizers for subcritical Trudinger-Moser functionals.
Let $\alpha<\lambda_1^{\mathbf{G}}$ be fixed. Then we have
\begin{lemma}\label{subcr}
 For any $0<\epsilon<4\pi\ell$, there exists some $u_\epsilon\in C^1(\Sigma\setminus
\{p_1,\cdots,p_L\},g_0)\cap C^0(\Sigma,g_0)\cap\mathscr{H}_{\mathbf{G}}$ with
$\|u_\epsilon\|_{1,\alpha}=1$ satisfying
\be\label{1-1}\int_\Sigma e^{(4\pi\ell-\epsilon)u_\epsilon^2}dv_g=\sup_{u\in\mathscr{H}_{\mathbf{G}},\|u\|_{1,\alpha}\leq 1}\int_\Sigma
e^{(4\pi\ell-\epsilon)u^2}dv_g.\ee
Moreover $u_\epsilon$ satisfies the Euler-Lagrange equation
\be\label{E-L-1}\le\{
  \begin{array}{lll}
  \triangle_g u_\epsilon-\alpha u_\e=
  \f{1}{\lambda_\epsilon}u_\epsilon e^{(4\pi \ell-\epsilon)u_\epsilon^2}-\f{\mu_\epsilon}{\lambda_\epsilon} \quad \mathrm{in}\,{\Si}, \quad \\[1.5ex]
    \int_{\Si}u_\epsilon dv_{g}=0,\\[1.5ex]
  \lambda_\epsilon=\int_\Si u_\epsilon^2 e^{(4\pi \ell-\epsilon)u_\epsilon^2}dv_g,\\[1.5ex]
  \mu_\epsilon=\f{1}{\mathrm{Vol}_g(\Si)}\int_\Si u_\epsilon e^{(4\pi \ell-\epsilon)u_\epsilon^2}dv_g,
  \end{array}
  \ri.\ee
  where $\Delta_g$ is the Laplace-Beltrami operator on $(\Sigma,g)$.
\end{lemma}
\proof Fix $\alpha<4\pi\ell$ and $0<\epsilon<4\pi\ell$.
Take a maximizing function sequence $(u_j)\subset \mathscr{H}_{\mathbf{G}}$ verifying that
$\|u_j\|_{1,\alpha}\leq 1$, and that as $j\ra\infty$,
$$
\ \int_{\Sigma} e^{(4\pi \ell-\epsilon)u_j^2}dv_g\ra \sup_{u\in \mathscr{H}_{\mathbf{G}},\,\|u\|_{1,\alpha}\leq 1}\int_\Si
e^{(4\pi \ell-\epsilon)u^2}dv_g.
$$
Clearly $(u_j)$ is bounded in $W^{1,2}(\Si,g)$. With no loss of generality we assume
$u_j$ converges to $u_\epsilon$  weakly  in
  $W^{1,2}(\Si,g)$,   strongly in  $L^s(\Si,g)$ for any $s>1$, and
  almost everywhere in $\Si$.
  This implies $u_\epsilon\in \mathscr{H}_{\mathbf{G}}$ and $\|u_\epsilon\|_{1,\alpha}\leq 1$. By Lemma \ref{b-2}, we have
  that  $e^{(4\pi \ell-\epsilon)u_j^2}$ converges to $e^{(4\pi \ell-\epsilon)u_\epsilon^2}$
  in $L^1(\Sigma,g)$ as $j\ra\infty$. Thus (\ref{1-1}) holds. It is easy to see that $\|u_\epsilon\|_{1,\alpha}=1$.

 By a simple calculation, $u_\epsilon$ is a distributional solution of the Euler-Lagrange equation (\ref{E-L-1}). In view of $g=\rho g_0$,
  applying elliptic estimates to (\ref{E-L-1}), we conclude  $u_\epsilon\in C^1(\Sigma\setminus\{p_1,\cdots,p_L\},g_0)\cap
  C^0(\Sigma,g_0)$. $\hfill\Box$\\

  Using the same argument as \cite[\,page\,3184]{Yang-JDE2015}, we get
  \be\label{lam-0}\liminf_{\epsilon\ra 0}\lambda_\epsilon>0,\quad |\mu_\epsilon|/\lambda_\epsilon\leq C.\ee
  \subsection{Blow-up analysis}
   Since $u_\epsilon$ is bounded in $W^{1,2}(\Sigma,g)$, we assume $u_\epsilon$ converges to some $u^\ast$ weakly in $W^{1,2}(\Sigma,g)$,
   strongly in $L^s(\Sigma,g)$ for any $s>1$, and almost everywhere in $\Sigma$. Obviously $\|u^\ast\|_{1,\alpha}\leq 1$.
   If $u_\epsilon$ is uniformly bounded, then by the Lebesgue dominated convergence theorem,
   \be\label{extremal}\int_\Sigma e^{4\pi\ell {u^\ast}^2}dv_g=\lim_{\epsilon\ra 0}\int_\Sigma e^{(4\pi\ell-\epsilon)u_\epsilon^2}dv_g=
   \sup_{u\in\mathscr{H}_{\mathbf{G}},\|u\|_{1,\alpha}\leq 1}\int_\Sigma e^{4\pi\ell u^2}dv_g.
   \ee
   Thus $u^\ast$ is the desired maximizer. In the following we assume $\max_\Sigma|u_\epsilon|\ra \infty$ as $\epsilon\ra 0$. Since
   $-u_\epsilon$ still satisfies (\ref{1-1}) and (\ref{E-L-1}), we assume with no loss of generality
   \be\label{bl-u}c_\epsilon=u_\epsilon(x_\epsilon)=\max_{\Sigma}|u_\epsilon|\ra\infty\ee
   and
   \be\label{x0}x_\epsilon\ra x_0\in\Sigma\ee
   as $\epsilon\ra 0$. To begin with, we have

\begin{lemma}\label{zero}
$u_\epsilon$ converges to $0$ weakly in $W^{1,2}(\Sigma,g)$, strongly in $L^s(\Sigma,g)$ for any $s>1$, and
almost everywhere in $\Sigma$.
\end{lemma}
\proof Since $u_\epsilon$ is bounded in $W^{1,2}(\Sigma,g)$, we assume $u_\epsilon$ converges to $u_0$ weakly in $W^{1,2}(\Sigma,g)$, strongly in $L^s(\Sigma,g)$ for any $s>1$, and
almost everywhere in $\Sigma$. Suppose $u_0\not\equiv 0$. Then
 $$\|u_\e-u_0\|_{1,\alpha}^2=\|u_\e\|_{1,\alpha}^2-\|u_0\|_{1,\alpha}^2+o_\e(1)\leq 1-\f{1}{2}\|u_0\|_{1,\alpha}^2$$
 for sufficiently small $\epsilon>0$. Using  the Young inequality, the H\"older inequality and  Lemma \ref{b-2}, we have that
 $e^{(4\pi\ell-\epsilon)u_\epsilon^2}$ is bounded in $L^q(\Sigma,g)$ for some $q>1$. Noting (\ref{lam-0}) and applying elliptic estimate
 to (\ref{E-L-1}), we obtain $u_\epsilon$ is uniformly bounded. This contradicts (\ref{bl-u}). Hence $u_0\equiv 0$.
 $\hfill\Box$\\

 Recalling the definitions of $I(x)$, $\beta(x)$ and $\ell$, namely (\ref{number})-(\ref{M0}), under the assumptions (\ref{bl-u}) and (\ref{x0}), we obtain the following energy concentration phenomenon. From now on, we write $ I_0=I(x_0)$ and $ \beta_0=\beta(x_0)$ for short, where $x_0$ is given as in (\ref{x0}).

\begin{lemma}\label{blow} $(i)$ $\lim_{r\ra 0}\lim_{\epsilon\ra 0}\int_{B_{g_0,r}(x_0)}|\nabla_{g_0} u_\e|^2dv_{g_0}=1/{I_0}$, where
$B_{g_0,r}(x_0)$ denotes the geodesic ball centered at $x_0$ with radius $r$ with respect to the metric $g_0$;
$(ii)$ $I_0(1+\beta_0)=\ell$.
\end{lemma}
\proof
We first prove the assertion $(i)$. With no loss of generality, we assume $\sigma_1(x_0),\cdots,\sigma_{I_0}(x_0)$ are all distinct points in $\mathbf{G}(x_0)$.
Choose some $r_0>0$ such that $B_{g_0,r_0}(\sigma_j(x_0))\cap B_{g_0,r_0}(\sigma_i(x_0))=\varnothing$ for every $1\leq i<j\leq I_0$.
Since $\int_\Sigma|\nabla_{g_0}u_\epsilon|^2dv_{g_0}=\int_\Si |\nabla_g u_\epsilon|^2 dv_g=1+o_\epsilon(1)$ and $ B_{g_0,r_0}(\sigma_k(x_0))=\sigma_k( B_{g_0,r_0}(x_0))$ for $k=1,\cdots,I_0$,
we have
\be\label{dirac1}
\int_{B_{g_0,r_0}(x_0)}|\nabla_{g_0} u_\e|^2dv_{g_0}\leq \f{1}{I_0}+o_\epsilon(1).
\ee
Suppose $(i)$ does not hold. There would exists a constant $\nu_0>0$ and $0<r_1<r_0$ such that
\be\label{dirac-2}
\int_{B_{g_0,r_1}(x_0)}|\nabla_{g_0} u_\e|^2dv_{g_0}\leq \f{1}{I_0}-\nu_0
\ee
for all sufficiently small $\epsilon>0$. Since
$\ell\leq \min\{I_0,I_0(1+\beta_0)\}\leq I_0$,
one finds a $p>1$ such that $e^{4\pi\ell u_\epsilon^2}$ is bounded in $L^{p}(B_{g_0,r_1/2}(x_0))$. In view of
(\ref{lam-0}) and Lemma \ref{zero}, one has by applying elliptic estimates to (\ref{E-L-1}) that
$u_\epsilon$ is bounded in $L^\infty(B_{g_0,r_1/4}(x_0))$, which contradicts the assumption (\ref{bl-u}).
This confirms $(i)$.

$(ii)$ Suppose not. Obviously $\ell<I_0(1+\beta_0)$. By (\ref{assumption}) and (\ref{betax}), we have
$\beta_0\leq 0$. This together with
 $(i)$ and an inequality of Adimurthi-Sandeep \cite{A-S} implies that
 there exist $r_0>0$,  $p>1$ and $C>0$ satisfying
 $$\int_{B_{g_0,r_0}(x_0)}e^{4\pi \ell p u_\epsilon^2}\rho dv_{g_0}\leq C.$$
 Applying elliptic estimates to (\ref{E-L-1}), we conclude that $u_\epsilon$ is bounded in
 $L^\infty(B_{g_0,r_0/2}(x_0))$, contradicting the assumption (\ref{bl-u}). Therefore $(ii)$ holds.
 $\hfill\Box$\\

 Set
\be\label{re}
r_\e=\sqrt{\lambda_\e}c_\e^{-1}e^{-(2\pi \ell-\e/2)c_\e^2}.
\ee
Using the same argument as that of derivation of (\cite{preprint}, the equation (42)), we have for any $0<a<4\pi\ell$,
\be\label{t-0}r_\epsilon^2c_\epsilon^2e^{(4\pi\ell-\epsilon-a)c_\epsilon^2}=o_\epsilon(1).\ee
In particular, $r_\epsilon\ra 0$ as $\epsilon\ra 0$. And it follows from (\ref{t-0}) that
\be\label{eta}r_\e^2 c_\e^q\ra 0,\quad \forall q>1.\ee

Keep in mind $g$ and $g_0$ satisfy (\ref{rho}), (\ref{h1}), and (\ref{isorho}). For any $1\leq k\leq N$, we  take an isothermal coordinate system $(U_{\sigma_k(x_0)},\psi_k;\{y_1,y_2\})$ near $\sigma_k(x_0)$  such that $\psi_k:U_{\sigma_k(x_0)}\ra
\Omega\subset\mathbb{R}^2$ is a homomorphism, $\psi_k(\sigma_k(x_0))=0$, and
\be \label{isocoordinate}
g_0=e^{2f_k}(dy_1^2+dy_2^2),
\ee
where $f_k\in C^1(\Omega,\mathbb{R})$ satisfies $f_k(0)=0$. If $g$ has a conical singularity of the order $\beta_0$ at $x_0$,
then in this coordinate system, $g$ can be represented by
\be\label{isosincoo}
g=V_ke^{2f_k}|y|^{2\beta_0} (dy_1^2 +dy_2^2 ),
 \ee
where $V_k\in C^0 (\Omega,\mathbb{R})$. It follows from  (\ref{h1}) and (\ref{isorho}) that
         \be\label{V0}
         V_k(0)=\lim_{d_{g_0}(x,x_0)\ra 0}\f{\rho(x)}{d_{g_0}(x,\sigma_k(x_0))^{2\beta_0}}=V_0,
         \ee
         where $V_0$ is a positive constant independent of $k$.
In particular, if  $\beta_0=0$, with no loss of generality, one can take $V_k(y)\equiv1$, and (\ref{isosincoo}) reduces to (\ref{isocoordinate}). Writing $\widetilde{x}_\epsilon=\psi_k^{-1}(x_\epsilon)$, we have the following:
\begin{lemma}\label{xerec}
If $\beta_0<0$, then $|\widetilde{x}_\e|^{1+\beta_0}/r_\e$ is uniformly bounded.
\end{lemma}
\proof For otherwise, up to a subsequence, we have
\be\label{case1}
|\widetilde{x}_\e|^{1+\beta_0}/r_\e\ra \infty.
\ee
For $y\in\Omega_{1,\epsilon}:=\{y\in \mathbb{R}^2:\widetilde{x}_\epsilon+r_\epsilon|\widetilde{x}_\e|^{-\beta_0} y\in\Omega\}$, we denote
 $$w_\epsilon(y)=c_\epsilon^{-1}({u}_\epsilon\circ\psi_k^{-1})(\widetilde{x}_\epsilon+r_\epsilon|\widetilde{x}_\e|^{-\beta_0} y),\quad
   v_\epsilon(y)=c_\epsilon \le(({u}_\epsilon\circ\psi_k^{-1})(\widetilde{x}_\epsilon+r_\epsilon|\widetilde{x}_\e|^{-\beta_0} y)-c_\e\ri).$$
  By (\ref{E-L-1}),  we calculate on $\Omega_{1,\epsilon}$,
  $$-\Delta_{\mathbb{R}^2} w_\epsilon=V_k(\widetilde{x}_\epsilon+r_\epsilon
   y)e^{2f_k(\widetilde{x}_\epsilon+r_\epsilon
   y)}|\widetilde{x}_\epsilon+r_\epsilon
   y|^{2\beta_0}|\widetilde{x}_\e|^{-2\beta_0}(\alpha r_\e^2w_\e+c_\epsilon^{-2}w_\epsilon e^{(4\pi \ell-\epsilon)c_\e^2(w_\e^2-1)}-c_\e^{-1}r_\e^2\mu_\e\lambda_\e^{-1}),$$
   $$-\Delta_{\mathbb{R}^2}v_\epsilon=V_k(\widetilde{x}_\epsilon+r_\epsilon
  y)e^{2f_k(\widetilde{x}_\epsilon+r_\epsilon
   y)}|\widetilde{x}_\epsilon+r_\epsilon
   y|^{2\beta_0}|\widetilde{x}_\e|^{-2\beta_0}(\alpha c_\e^2 r_\e^2w_\e+w_\epsilon e^{(4\pi \ell-\epsilon)(1+w_\epsilon)v_\epsilon}
  -c_\e r_\e^2\mu_\e\lambda_\e^{-1}),
  $$
   where $\Delta_{\mathbb{R}^2}$ stands for the standard Laplacian operator on $\mathbb{R}^2$.
   It follows from (\ref{case1}) that both $e^{2f(\tilde{x}_\epsilon+r_\epsilon
   y)}$ and $|\tilde{x}_\epsilon+r_\epsilon
   y|^{2\beta_0}\cdot|\tilde{x}_\e|^{-2\beta_0}$ are $1+o_\epsilon(1)$ in $\mathbb{B}_R$ for any fixed $R>0$. Combining (\ref{lam-0}) with (\ref{eta}), we have $c_\e r_\e^2\mu_\e\lambda_\e^{-1}$ is $o_\e(1)$. Applying elliptic estimates to the above two
    equations, we obtain
  \be\label{loc-w1}w_\epsilon\ra1 \quad {\rm in}\quad C^1_{\rm loc}(\mathbb{R}^2)\ee
  and $v_\e\ra v_0$ in $C^1_{\mathrm{loc}}(\mathbb{R}^2)$, where $v_0$ satisfies
\be\label{v0}
\le\{
     \begin{array}{llll}
     &-\Delta_{\mathbb{R}^2} v_0=V_0e^{8\pi \ell v_0}\qquad \textrm{in}\,\mathbb{R}^2 \\[1.5ex]
     &v_0(0)=0=\sup_{\mathbb{R}^2}v_0
     \end{array}
     \ri.
\ee
in the is distributional sense.
Moreover, one easily estimates
\be\label{v0int1}
\int_{\mathbb{R}^2}e^{8\pi\ell v_0}\mathrm{d}y\leq \f{1}{I_0V_0}.
\ee
In view of (\ref{v0}) and (\ref{v0int1}), we have by a classification theorem of Chen-Li \cite{Chen-Li},
    $$v_0(y)=-\f{1}{4\pi \ell}\log(1+\pi \ell V_0|y|^2).$$
    It then follows that
    \be\label{v1m0}
    \int_{\mathbb{R}^2}e^{8\pi v_0}\mathrm{d}y=\f{1}{\ell V_0}.
    \ee
    Since $\beta_0<0$, it follows from $(ii)$ of Lemma \ref{blow} that $\ell<I_0$. As a consequence,
    there is a contradiction between
(\ref{v1m0}) and (\ref{v0int1}). This ends the proof of the lemma.
 $\hfill\Box$\\

 We now define two sequences of functions
 \be\label{psie,varphie}
\psi_\e(y)=c_\e^{-1}\widetilde{u}_\e (\widetilde{x_\e} +r_\e^{1/(1+\beta_0)}y),\quad \varphi_\e(y)=c_\e \le(\widetilde{u}_\e
(\widetilde{x}_\e +r_\e^{1/(1+\beta_0)}y)-c_\e\ri)
\ee
for  $y\in\Omega_{2,\epsilon}:=\le\{y\in \mathbb{R}^2:\widetilde{x}_\epsilon+r_\e^{1/(1+\beta_0)}y\in\Omega\ri\}$. Then there holds
the following:
\begin{lemma}\label{blow-upanal} If $\beta_0<0$, then
$(i)$ $\psi_\epsilon\ra 1$ in $C^0_{\mathrm{loc}}({\mathbb{R}^2})\cap W_{\rm loc}^{1,2}(\mathbb{R}^2)$;
$(ii)$ $\varphi_\e\ra \varphi $ in $C^0_{\mathrm{loc}}({\mathbb{R}^2})\cap W_{\rm loc}^{1,2}(\mathbb{R}^2)$, where
\be\label{b-b-0}\varphi(y)=-\f{1}{4\pi \ell}\log\le(1+\f{\pi I_0V_0}{1+\beta_0}|y|^{2(1+\beta_0)}\ri).\ee
\end{lemma}
\proof
  By (\ref{E-L-1}) and (\ref{re}),  we have on $\Omega_{2,\epsilon}$,
  \bea
  {\nonumber}
 -\Delta_{\mathbb{R}^2} \psi_\epsilon&=&V_k(\widetilde{x}_\epsilon+r_\e^{1/(1+\beta_0)}y)e^{2f_k(\widetilde{x}_\epsilon+r_\e^{1/(1+\beta_0)}y)}|y+r_\e^{-1/(1+\beta_0)}
 \widetilde{x}_\epsilon|^{2\beta_0}(\alpha r_\e^2\psi\\
 &&\quad+c_\epsilon^{-2}\psi_\epsilon e^{(4\pi \ell-\epsilon)c_\epsilon^2(\psi_\e^2-1)}-c_\e^{-1}r_\e^{2}\mu_\e\lambda_\e^{-1}), \label{psi1}\\
  {\nonumber}-\Delta_{\mathbb{R}^2}\varphi_\epsilon&=&V_k(\widetilde{x}_\epsilon+r_\e^{1/(1+\beta_0)}y)
  e^{2f_k(\widetilde{x}_\epsilon+r_\e^{1/(1+\beta_0)}y)}|y+r_\e^{-1/(1+\beta_0)}\widetilde{x}_\epsilon|^{2\beta_0}(\alpha c_\e^2 r_\e^2\psi_\e\\
  \label{phi1}&&\quad+\psi_\epsilon e^{(4\pi \ell-\epsilon)(1+\psi_\epsilon)\varphi_\epsilon}-c_\e r_\e^{2}\mu_\e\lambda_\e^{-1}).
   \eea
  In view of Lemma \ref{xerec}, $r_\e^{-1/(1+\beta_0)}\widetilde{x}_\epsilon$ is a bounded sequence of points. We may assume with no loss of generality that
   $r_\e^{-1/(1+\beta_0)}\widetilde{x}_\epsilon\ra p\in\mathbb{R}^2$ as $\epsilon\ra 0$. Note that $\beta>-1$. Applying elliptic estimates to (\ref{psi1}),
    we obtain $\psi_\e\ra \psi$ in $C^0_{\mathrm{loc}}({\mathbb{R}^2})\cap W_{\rm loc}^{1,2}(\mathbb{R}^2)$, where $\psi$ is a distributional harmonic function. Then the Liouville theorem leads to $\psi\equiv 1$.
    Further application of elliptic estimates on (\ref{phi1}) implies that
\be\label{loc-phi1}
\varphi_\epsilon\ra \varphi  \quad \text{in}\quad C^0_{\mathrm{loc}}({\mathbb{R}^2})\cap W_{\rm loc}^{1,2}(\mathbb{R}^2),
\ee
where $\varphi$ is a distributional solution of
\be\label{varphi0}
\le\{
     \begin{array}{llll}
     &-\Delta_{\mathbb{R}^2} \varphi=|y+p|^{2\beta_0}V_0e^{8\pi I_0(1+\beta_0)\varphi}\quad {\rm in}\quad\mathbb{R}^2\\[1.5ex]
     &\varphi(0)=0=\max_{\mathbb{R}^2}\varphi.
     \end{array}
     \ri.
\ee
One can easily derive
\be\label{energ}
\int_{\mathbb{R}^2}V_0|y+p|^{2\beta_0}e^{8\pi \ell\varphi}dy\leq1.
\ee
In view of (\ref{varphi0}) and (\ref{energ}), a classification theorem of Chen-Li
\cite{CL-sin} suggests the representation:
\be\label{v0solu}
\varphi(y)=-\f{1}{4\pi \ell}\log\le(1+\f{\pi I_0V_0}{1+\beta_0}|y+p|^{2(1+\beta_0)}\ri).
\ee
Since $\varphi(0)=0$, we have $p=0$. By a straightforward calculation,
\be\label{v0soluint}
\int_{\mathbb{R}^2}V_0|y|^{2\beta_0}e^{8\pi \ell\varphi(y)}dy=\f{1}{I_0},
\ee
as desired.
 $\hfill\Box$\\

In the case $\beta_0=0$, we have an analog of  Lemma \ref{blow-upanal}, namely
\begin{lemma}\label{beta=0}
Let $\psi_\epsilon$ and $\varphi_\epsilon$ be defined as in (\ref{psie,varphie}). If $\beta_0=0$, then
$\psi_\epsilon\ra 1$ and
$\varphi_\e\ra \varphi $ in $C^1_{\mathrm{loc}}({\mathbb{R}^2})$, where
$\label{b-b-0}\varphi(y)=-\f{1}{4\pi I_0}\log\le(1+{\pi I_0\rho(x_0)}|y|^{2}\ri)$, $\rho$ is given as in
(\ref{rho}) and (\ref{h1}).
\end{lemma}

\proof Noting that $\beta_0=0$, we have by applying elliptic estimates to (\ref{psi1}) and (\ref{phi1}) that
$\psi_\epsilon\ra 1$ and
$\varphi_\e\ra \varphi $ in $C^1_{\mathrm{loc}}({\mathbb{R}^2})$, where $\varphi$ satisfies
$$
\le\{
     \begin{array}{llll}
     -\Delta_{\mathbb{R}^2} \varphi=\rho(x_0)e^{8\pi I_0\varphi}\quad {\rm in}\quad\mathbb{R}^2\\[1.5ex]
     \varphi(0)=0=\max_{\mathbb{R}^2}\varphi\\[1.5ex]
     \int_{\mathbb{R}^2}\rho(x_0)e^{8\pi I_0\varphi}dy\leq1.
     \end{array}
     \ri.
$$
Then a result of Chen-Li \cite{Chen-Li} leads to $\varphi(y)=-\f{1}{4\pi I_0}\log\le(1+{\pi I_0\rho(x_0)}|y|^{2}\ri)$. As a consequence,
\be\label{b-0}\int_{\mathbb{R}^2}\rho(x_0)e^{8\pi I_0\varphi(y)}dy=\f{1}{I_0},\ee
which is an analog of (\ref{v0soluint}). $\hfill\Box$\\

By (\ref{re}), Lemmas  \ref{blow-upanal} and \ref{beta=0}, we have for any fixed $R>0$,
\bna
{\nonumber}\int_{\mathbb{B}_R(0)}V_0|y|^{2\beta_0}e^{8\pi \ell\varphi}dy&=&\lim_{\e\ra0}
\int_{\mathbb{B}_R(0)}V_0|y|^{2\beta_0}e^{(4\pi \ell-\e)(1+\psi_\e)\varphi_\e}dy\\
{\nonumber}&=&\lim_{\e\ra0}\f{1}{\lambda_{\e}} \int_{\mathbb{B}_{R r_\e^{1/(1+\beta_0)}}(\widetilde{x}_\epsilon)}V_0e^{2f_k}|y|^{2\beta_0}\widetilde{u}_\e^2 e^{(4\pi \ell-\e)\widetilde{u}_{\e}^2}dy\\
&=&\lim_{\e\ra0}\f{1}{\lambda_{\e}}\int_{\psi_k^{-1}(\mathbb{B}_{R r_\e^{1/(1+\beta_0)}}(\widetilde{x}_\epsilon))}u_\e^2 e^{(4\pi \ell-\e)u^2_{\e}}dv_g,
\ena
where $V_0=\rho(x_0)$ and $\ell=I_0$ if $\beta_0=0$. This together with (\ref{v0soluint}) and (\ref{b-0})
implies that
\be\label{bubble}
\lim_{R\ra\infty}\lim_{\e\ra 0}\f{1}{\lambda_{\e}}\int_{\psi_k^{-1}(\mathbb{B}_{R r_\e^{1/(1+\beta_0)}}(\widetilde{x}_\epsilon))}u_\e^2 e^{(4\pi \ell-\e)u^2_{\e}}dv_g=\f{1}{I_0}.
\ee
Noting that
\bna
\lambda_\epsilon
=\int_{\cup_{k=1}^{I_0}\psi_k^{-1}(\mathbb{B}_{R r_\e^{1/(1+\beta_0)}}(\widetilde{x}_\epsilon))}u_\epsilon ^2e^{(4\pi \ell-\e)u_\epsilon^2}dv_g+\int_{\Si\setminus \cup _{k=1}^{I_0}\psi_k^{-1}(\mathbb{B}_{Rr_\e^{1/(1+\beta_0)}}(\widetilde{x}_{\epsilon}))}u_\epsilon ^2e^{(4\pi \ell-\e)u_\epsilon^2}dv_g,
\ena
we conclude from (\ref{bubble}) that
\be\label{outbubble}
\lim_{R\ra\infty}\lim_{\e\ra0}\f{1}{\lambda_\e}\int_{\Si\setminus \cup _{k=1}^{I_0}\psi_k^{-1}(\mathbb{B}_{Rr_\e^{1/(1+\beta_0)}}(\widetilde{x}_{\epsilon}))}u_\epsilon ^2e^{(4\pi \ell-\e)u_\epsilon^2}dv_g=0.
\ee

Similar to \cite{Lijpde}, we define $u_{\e,\gamma}=\min \{u_\e,\gamma c_\e\}$ for any $0<\gamma<1$, and have
\begin{lemma} For any $0<\gamma<1$, there holds
$$\label{gama}
\lim_{\e\ra0}\int_{\Si}|\nabla_g u_{\e,\gamma} |^2d v_g=\gamma.
$$
\end{lemma}
 \proof For fixed $R>0$ and  sufficiently small $\e$, in view of (\ref{E-L-1}), we have by using integration by parts,
  (\ref{bubble}) and (\ref{outbubble}) that
\bna
\int_{\Si}|\nabla_g u_{\e,\gamma} |^2d v_g&=&\int_{\Si}\nabla_g u_{\e,\gamma} \nabla_g u_\e dv_g\\
&=&\lambda_\e^{-1}\int_{\cup_{k=1}^{I_0}\psi_k^{-1}(\mathbb{B}_{R r_\e^{1/(1+\beta_0)}}(\widetilde{x}_\epsilon))}u_\epsilon u_{\e,\gamma}e^{(4\pi \ell-\e)u_\epsilon^2}dv_g\\
&&+\lambda_\e^{-1}\int_{\Si\setminus \cup _{k=1}^{I_0}\psi_k^{-1}(\mathbb{B}_{Rr_\e^{1/(1+\beta_0)}}(\widetilde{x}_{\epsilon}))}u_\epsilon u_{\e,\gamma}e^{(4\pi \ell-\e)u_\epsilon^2}dv_g+o_\e(1)\\
&=& (1+o_\e(1))\gamma{\lambda_\e}^{-1}\int_{\cup_{k=1}^{I_0}\psi_k^{-1}(\mathbb{B}_{R r_\e^{1/(1+\beta_0)}}(\widetilde{x}_\epsilon))}u_\epsilon^2 e^{(4\pi \ell-\e)u_\epsilon^2}dv_g+o(1)\\
&=&\gamma+o(1),
\ena
where $o(1)\ra 0$ as $\epsilon\ra 0$ first, and then $R\ra\infty$.
The lemma follows immediately.{$\hfill\Box$}

\begin{lemma}\label{tends-0}
There holds $c_\epsilon/\lambda_\epsilon\ra 0$ as $\epsilon\ra 0$.
\end{lemma}
\proof For any fixed $0<\gamma<1$,
\bea
{\nonumber}\int_{\Si}e^{( 4\pi \ell-\e) u_\epsilon^2}dv_g&=&\int_{u_\e\leq \gamma c_\e}e^{( 4\pi \ell-\e)u_\epsilon^2}dv_g+\int_{u_\e> \gamma c_\e}e^{( 4\pi \ell-\e)u_\epsilon^2}dv_g\\[1.2ex]
\label{gamma-3}&\leq&\int_{\Si}e^{( 4\pi \ell-\e) u_{\epsilon,\gamma}^2}dv_g+\dfrac{\lambda_\e}{\gamma^2c_\e^2}.
\eea
By Lemmas \ref{zero} and \ref{gama}, we conclude
$$\int_\Sigma e^{( 4\pi \ell-\e) u_{\epsilon,\gamma}^2}dv_g=\mathrm{Vol}_g(\Si)+o_\epsilon(1).$$
Passing to the limit $\epsilon\ra 0$ first, and then $\gamma\ra 1$ in (\ref{gamma-3}), we have
\be\label{gama2}
\sup_{u\in \mathscr{H}_{\mathbf{G}},\,\|u\|_{1,\alpha}\leq 1}\int_{\Si}e^{4\pi \ell u^2}dv_g=\lim_{\epsilon\ra 0}
\int_\Sigma e^{4\pi\ell u_\epsilon^2}dv_g\leq \mathrm{Vol}_g(\Si)+\liminf_{\e\ra0}\f{\lambda_\e}{c_\e^2}.
\ee
Since
$$
\sup_{u\in \mathscr{H}_{\mathbf{G}},\,\|u\|_{1,\alpha}\leq 1}\int_{\Si}e^{4\pi \ell u^2}dv_g> \mathrm{Vol}_g(\Si),
$$
we have by (\ref{gama2}) that $\liminf_{\e\ra0}{\lambda_\e}/{c_\e^2}
>0$. In particular $c_\epsilon/\lambda_\epsilon\ra 0$ as $\epsilon\ra 0$. $\hfill\Box$\\

{Recall $\mathbf{G}(x_0)=\{\sigma_1(x_0),\cdots,\sigma_{I_0}(x_0)\}$, and  $\mathbf{S}=\{p_1,\cdots,p_L \}$.} The convergence of $c_\epsilon u_\epsilon$ is precisely described as follows.
 \begin{lemma}\label{converge-Green} For any $1<q<2$, we have
 $c_\epsilon u_\epsilon$ converges to $G_\alpha$ weakly in $W^{1,q}(\Sigma,g_0)$, strongly in $L^{2q/(2-q)}(\Sigma)$,
 and in ${ C^1(\Sigma\setminus\{\mathbf{G}(x_0)\cup\mathbf{S}\})}$, where $G_\alpha$ is a Green function satisfying
 \be\label{Green-01}\le\{\begin{array}{lll}
 \Delta_{g_0}G_\alpha-\alpha\rho G_\alpha=\f{1}{I_0}\sum_{i=1}^{I_0} \delta_{\sigma_i(x_0)}-\f{\rho}{{\rm Vol}_g(\Sigma)}\\[1.5ex]
 \int_\Sigma G_\alpha dv_g=0\\[1.5ex]
 G_\alpha(\sigma_i(x))=G_\alpha(x),\, x\in \Sigma\setminus \{\sigma_j(x_0)\}_{j=1}^{I_0},\, 1\leq i\leq I_0.
 \end{array}\ri.\ee
 \end{lemma}
 \proof In view of (\ref{E-L-1}), one has
 \be\label{ge1}
\le\{
     \begin{array}{llll}
\Delta_{g} (c_\e u_\epsilon)-\alpha ( c_\e u_\e)=
 f_\e-b_\e\,\,\,{\rm on}\,\,\, \Si\\[1.2ex]
  \int_{\Si}c_\e u_\e dv_{g}=0\\[1.2ex]
  f_\epsilon=\f{1}{\lambda_\epsilon}  c_\e u_\e e^{( 4\pi \ell-\e) u_\epsilon^2}\\[1.2ex]
  b_\epsilon=\f{ c_\e\mu_\e}{\lambda_\e}.
  \end{array}\ri.
\ee

Firstly we claim that
\be\label{weak}f_\epsilon dv_g\rightharpoonup \f{1}{I_0}\sum_{i=1}^{I_0} \delta_{\sigma_i(x_0)}\ee
weakly in the sense of measure, or equivalently, there holds
$$\int_\Sigma f_\epsilon\phi dv_g=\f{1}{I_0}\sum_{i=1}^{I_0}\phi(\sigma_i(x_0))+o_\epsilon(1),\quad\forall \phi\in C^0(\Sigma,g_0).$$
To see it, we estimate for any fixed $0<\gamma<1$ and $R>0$
\bea
{\nonumber}\int_{\Si}f_\e\phi dv_g&=&\int_{u_\e\leq \gamma c_\e}f_\e\phi dv_g+\int_{\{u_\e> \gamma c_\e\}\cap \cup_{k=1}^{I_0} \psi_k^{-1}(\mathbb{B}_{Rr_\e^{1/(1+\beta_0)}}(\widetilde{x}_\epsilon))}f_\e\phi dv_g\\[1.2ex]
{\nonumber}&\,&+\int_{\{u_\e> \gamma c_\e\}\backslash \cup_{k=1}^{I_0}\psi_k^{-1}(\mathbb{B}_{Rr_\e^{1/(1+\beta_0)}}(\widetilde{x}_\epsilon))}f_\e\phi dv_g\\[1.2ex]
\label{123}&:=& {\rm I}+{\rm II}+{\rm III}.
\eea
By Lemmas \ref{zero}, \ref{gama} and \ref{tends-0}, we have by the H\"older inequality
$${\rm I}=\f{c_\epsilon}{\lambda_\epsilon}\int_{u_\e\leq \gamma c_\e}u_\e e^{ (4\pi \ell-\e)u_{\e,\gamma}^2}\phi dv_g=o_\e(1).
$$
Note that $\cup_{k=1}^{I_0}\psi_k^{-1}(\mathbb{B}_{Rr_\e^{1/(1+\beta_0)}}(\widetilde{x}_\epsilon))\subset\{u_\e> \gamma c_\e\}$ for  sufficiently small $\e>0$. In view of Lemmas \ref{blow-upanal} and \ref{beta=0}, we calculate by using
(\ref{bubble}) and the mean value theorem for integrals,
\bna
{\rm II}&=&\int_{\cup_{k=1}^{I_0}\psi_k^{-1}(\mathbb{B}_{Rr_\e^{1/(1+\beta_0)}}(\widetilde{x}_\epsilon))}f_\e\phi dv_g\\[1.2ex]
&=&\sum_{k=1}^{I_0}\phi\le(\sigma_k(x_0)\ri)\le(1+o_\e(1)\ri)\int_{\psi_k^{-1}(\mathbb{B}_{Rr_\e^{1/(1+\beta_0)}}
(\widetilde{x}_\epsilon))}\f{c_\e}{\lambda_\e}
u_\e e^{ (4\pi
\ell-\e) u_\e^2} dv_g\\[1.2ex]
&=&\sum_{i=1}^{I_0}\phi(\sigma_i(x_0))\le(\f{1}{I_0}+o(1)\ri),
\ena
and
\bna
{\rm III}&\leq&\int_{\{u_\e> \gamma c_\e\}\setminus \cup_{k=1}^{I_0}\psi_k^{-1}(\mathbb{B}_{Rr_\e^{1/(1+\beta_0)}}(\widetilde{x}_\epsilon))}f_\e|\phi| dv_g\\[1.2ex]
&\leq&\f{\sup_\Sigma|\phi|}{\gamma}\int_{\{u_\e> \gamma c_\e\}\setminus \cup_{k=1}^{I_0}\psi_k^{-1}(\mathbb{B}_{Rr_\e^{1/(1+\beta_0)}}(\widetilde{x}_\epsilon))}\lambda_\e^{-1}u_\e^2 e^{( 4\pi \ell-\e) u_{\e}^2} dv_g \\[1.2ex]
&\leq&\f{\sup_\Sigma|\phi|}{\gamma}\le(1-\int_{\cup_{k=1}^{I_0}\psi_k^{-1}(\mathbb{B}_{Rr_\e^{1/(1+\beta_0)}}
(\widetilde{x}_\epsilon))}\lambda_\e^{-1}u_\e^2 e^{( 4\pi \ell-\e) u_{\e}^2} dv_g\ri)\\[1.2ex]
&=&o(1),
\ena
where $o(1)\ra 0$ as $\epsilon\ra 0$ first, and then $R\ra\infty$.
Inserting the estimates of ${\rm I}$-${\rm III}$ into (\ref{123}), we conclude our claim (\ref{weak}).

Secondly we calculate $b_\epsilon$ in (\ref{ge1}). Similar to the estimate of (\ref{123}), we have for any fixed $0<\gamma<1$,
\bna
\f{c_\epsilon}{\lambda_\epsilon}\int_{u_\epsilon\leq \gamma c_\epsilon} u_\epsilon e^{(4\pi\ell-\epsilon)u_\epsilon^2}dv_g=o_\epsilon(1)
\ena
and
\bna
\f{c_\epsilon}{\lambda_\epsilon}\int_{u_\epsilon> \gamma c_\epsilon}u_\epsilon e^{(4\pi\ell-\epsilon)u_\epsilon^2}dv_g
&=&\int_{\cup_{k=1}^{I_0}\psi_k^{-1}(\mathbb{B}_{Rr_\e^{1/(1+\beta_0)}}
(\widetilde{x}_\epsilon))}\f{1}{\lambda_\epsilon}u_\epsilon^2e^{(4\pi\ell-\epsilon)u_\epsilon^2}dv_g+o(1)\\
&=&1+o(1).
\ena
It then follows that
\be\label{b-epsilon}b_\epsilon=\f{1}{{\rm Vol}_g(\Sigma)}\f{c_\epsilon}{\lambda_\epsilon}\int_\Sigma u_\epsilon e^{(4\pi\ell-\epsilon)u_\epsilon^2}dv_g=\f{1}{{\rm Vol}_g(\Sigma)}+o_\epsilon(1).\ee

Thirdly we prove that $c_\epsilon u_\epsilon$ is bounded in $L^1(\Sigma,g)$. Suppose on the contrary
\be\label{Cont}\|c_\epsilon u_\epsilon\|_{L^1(\Sigma,g)}\ra \infty.\ee
Since for any fixed $0<\gamma<1$,
$$\int_\Sigma |f_\epsilon|dv_g=\int_{u_\epsilon\leq \gamma c_\epsilon}|f_\epsilon|dv_g+\int_{u_\epsilon>\gamma c_\epsilon}f_\epsilon dv_g,$$
we have that $f_\epsilon$ is bounded in $L^1(\Sigma,g)$ by using a similar argument of the estimate of (\ref{123}).
Obviously $b_\epsilon$ is a bounded sequence of numbers due to (\ref{b-epsilon}).
Define $w_\epsilon=c_\epsilon u_\epsilon/\|c_\epsilon u_\epsilon\|_{L^1(\Sigma,g)}$. Then (\ref{ge1}) gives
\be\label{ge-01}
\le\{
     \begin{array}{llll}
\Delta_{g_0} w_\epsilon=h_\epsilon:=\alpha\rho w_\e+ \rho\f{f_\e-b_\e}{\|c_\epsilon u_\epsilon\|_{L^1(\Sigma,g)}}\,\,\,{\rm on}\,\,\, \Si\\[1.2ex]
  \int_{\Si}w_\e dv_{g}=0\\[1.2ex]
  \|w_\epsilon\|_{L^1(\Sigma,g)}=1.
  \end{array}\ri.
\ee
Clearly we have got
\be\label{h-L1}\int_\Sigma|h_\epsilon|dv_{g_0}\leq C.\ee
 By the Green representation formula,
\be\label{Green-repres}w_\epsilon(x)-\f{1}{{\rm Vol}_{g_0}(\Sigma)}\int_\Sigma w_\epsilon dv_{g_0}=\int_\Sigma G(x,y)h_\epsilon(y)dv_{g_0,y},\ee
where $G(x,y)$ is the Green function for $\Delta_{g_0}$. In particular there exists a constant $C$ such that $|G(x,y)|\leq C{\rm dist}_{g_0}(x,y)$ and
$|\nabla_{g_0,x}G(x,y)|\leq C({\rm dist}_{g_0}(x,y))^{-1}$ for all $x,y\in\Sigma$. By (\ref{assumption}),
$\rho(x)$ has a positive
lower bound on $\Sigma$. As a consequence
\be\label{mean}\f{1}{{\rm Vol}_{g_0}(\Sigma)}\int_\Sigma |w_\epsilon| dv_{g_0}\leq C\int_\Sigma |w_\epsilon|\rho dv_{g_0}=C.\ee
Combining (\ref{h-L1}) and (\ref{Green-repres}), we obtain for any $1<q<2$,
$$\int_\Sigma|\nabla_{g_0}w_\epsilon|^qdv_{g_0}\leq C\int_\Sigma |h_\epsilon|dv_{g_0}\leq C.$$
While (\ref{Green-repres}) and (\ref{mean}) imply that for any $q>1$, there holds $\|w_\e\|_{L^q(\Sigma,g_0)}\leq C$.
Therefore $w_\epsilon$ is bounded in $W^{1,q}(\Sigma,g_0)$ for any $1<q<2$. The Sobolev embedding theorem leads to
$w_\e$ converges to $w$ weakly in $W^{1,q}(\Sigma,g_0)$, strongly in $L^r(\Sigma,g_0)$ for any $r<{2q}/{(2-q)}$, and
almost everywhere in $\Sigma$. Clearly $w$ satisfies
$$\le\{\begin{array}{lll}\Delta_{g_0}w=\alpha\rho w\quad{\rm in}\quad \Sigma\\[1.2ex]
\int_\Sigma w\rho dv_{g_0}=0.\end{array}\ri.$$
Since ${\rho\in L^r(\Sigma\setminus\bigcup_{i=1}^{L} B_{g_0,\d}(p_i),g_0)}$ { for any small $\d>0$ and some $r>1$} , we have $w\in { C^1(\Sigma\setminus\mathbf{S},g_0)}$ and $u\in\mathscr{H}_{\bm G}$ by using elliptic estimates.
Then integration by parts gives
$$\int_\Sigma|\nabla_{g}w|^2dv_{g}=\alpha\int_{\Sigma}w^2dv_g,$$
which leads to $w\equiv 0$ due to $\alpha<\lambda_1^{\bm G}$. This contradicts
$\|w\|_{L^1(\Sigma,g)}=\lim_{\epsilon\ra 0}\|w_\epsilon\|_{L^1(\Sigma,g)}=1$. Therefore $c_\epsilon u_\epsilon$
is bounded in $L^1(\Sigma,g)$.

Fourthly we analyze the convergence of $c_\epsilon u_\epsilon$. Rewrite (\ref{ge1}) as
\be\label{ge-30}
\le\{
     \begin{array}{llll}
\Delta_{g_0} (c_\e u_\epsilon)=\xi_\epsilon:=\alpha \rho c_\e u_\e+
 \rho(f_\e-b_\e)\,\,\,{\rm on}\,\,\, \Si\\[1.2ex]
  \int_{\Si}c_\e u_\e\rho dv_{g_0}=0.
  \end{array}\ri.
\ee
Now since $\xi_\epsilon$ is bounded in $L^1(\Sigma,g_0)$, we conclude that $c_\epsilon u_\epsilon$ is bounded in
$W^{1,q}(\Sigma,g_0)$ for any $1<q<2$ similar to $w_\epsilon$. Hence $c_\epsilon u_\epsilon$ converges to some $G_\alpha$ weakly in $W^{1,q}(\Sigma,g_0)$, strongly in $L^r(\Sigma,g_0)$ for any $r<2q/(2-q)$, and almost everywhere
in $\Sigma$. In view of (\ref{weak}) and (\ref{b-epsilon}), $G_\alpha$ satisfies (\ref{Green-01}) in the distributional sense.
Applying elliptic estimates to (\ref{ge-30}), we have that $c_\e u_\e$ converges to $G_\alpha$ in ${C^1(\Sigma\setminus\{\mathbf{G}(x_0)\cup\mathbf{S}\})}$.
This completes the proof of the lemma. $\hfill\Box$

\subsection{Upper bound estimate}
Recall  the isothermal coordinate system $(U_{\sigma_k(x_0)},\psi_k)$ near $\sigma_k(x_0)$ (here we only take $k$ from 1 to $I_0$) given as in (\ref{isocoordinate}). Set
  $$\label{r0}r_0=\f{1}{4}\min_{1\leq i<j\leq I_0}d_{g_0}(\sigma_i(x_0),\sigma_j(x_0)).$$
  For $\delta<r_0$ with $B_{g_0,2\d}(x_0)\subset U_{x_0}$, there exists two positive constants $c_1(\d)$ and $c_2(\d)$ such that $B_{g_0,(1-c_1(\d))\delta}(\sigma_k(x_0))\subset\psi_k^{-1}(\mathbb{B}_{\d})\subset B_{g_0,(1+c_2(\d))\delta}(\sigma_k(x_0)).$ Moreover, both $c_1(\d)$ and $c_2(\d)$ converge to 0 as $\d\ra0$. Hence, on $B_{g_0,2\delta}(\sigma_k(x_0))$, by using isothermal coordinate system $(U_k,\psi_k)$, (\ref{Green-01}) can be rewritten as the equation $G_\alpha\circ\psi^{-1}_k$ satisfies on $\psi_k(B_{g_0,2\delta}(\sigma_k(x_0)))$. By using elliptic estimates to that equation, we obtain:
\be\label{Glocphi}
G_\alpha\circ\psi^{-1}_k=-\f{1}{2\pi I_0}\log|y|+A_0+{{\Psi}_k}(y),
\ee
where ${{\Psi}_k}\in C^1(\mathbb{B}_{\f{5}{3}\d})$ satisfies ${{\Psi}_k}(0)=0$ for small $\d$, and $A_0$ is a constant defined by
\be
\label{a-0}A_0=\lim_{y\ra 0}\left(G_\alpha\circ\psi^{-1}_k(y)+\f{1}{2\pi I_0}\log |y|\right)
=\lim_{x\ra x_0}\left(G_\alpha(x)+\f{1}{2\pi I_0}\log d_{g_0}(x,x_0)\right).
\ee
By (\ref{Glocphi}), $G_\alpha$ near $x_0$ can be locally presented by
 \be\label{Gloc0}
G_\alpha(x)=-\f{1}{2\pi I_0}{\log d_{g_0}(x,x_0)}+A_0+\widetilde{{\Psi}}(x),
\ee
where $\widetilde{{\Psi}}\in C^1(B_{g_0,\f{3}{2}\d}(x_0))$ satisfies $\widetilde{{\Psi}}(x_0)=0.$ Furthermore, we obtain $G_\alpha$ near $\sigma_k(x_0)$ can be locally presented by
\be\label{Gloc}
G_\alpha(x)=-\f{1}{2\pi I_0}{\log d_{g_0}(x,\sigma_k(x_0))}+A_0+\widetilde{{\Psi}}(\sigma_k^{-1}(x)).
\ee
This conclusion is based on an observation that, for $x\in{B}_{g_0,\f{3}{2}\d}(\sigma_k(x_0))$, by (\ref{Gloc0})
\bna
{\nonumber}G_\alpha(x)+\f{1}{2\pi I_0}{\log d_{g_0}(x,\sigma_k(x_0))}-A_0&=&\big(G_\alpha(\sigma_k^{-1}(x))+\f{1}{2\pi I_0}{\log d_{g_0}(\sigma_k^{-1}(x),x_0)}-A_0\big)\\[1.2ex]
\label{Phisym}&=&\widetilde{{\Psi}}(\sigma_k^{-1}(x)).
\ena
{In the view of (\ref{E-L-1}), integration by parts leads to
\bna
{\nonumber}\int_{\Si\setminus \bigcup_{k=1}^{I_0} \psi_k^{-1}(\mathbb{B}_{\delta})}|\nabla_g u_\e|^2 d v_g&=&\int_{\Si\setminus \bigcup_{k=1}^{I_0} \psi_k^{-1}(\mathbb{B}_{\delta})}|\nabla_{g_0} u_\e|^2 d v_{g_0}\\[1.2ex]
{\nonumber}&=&- \sum_{k=1}^{I_0}\int_{ \partial\psi_k^{-1}(\mathbb{B}_{\delta})}u_\e\frac{\p u_\e}{\p n} ds_{g_0}+\int_{\Si\setminus \bigcup_{k=1}^{I_0} \psi_k^{-1}(\mathbb{B}_{\delta})}u_\e\triangle_{g_0} u_\e d v_{g_0}\\[1.2ex]
\label{gintd}&=&- \sum_{k=1}^{I_0}\int_{ \partial\psi_k^{-1}(\mathbb{B}_{\delta})}u_\e\frac{\p u_\e}{\p n} ds_{g_0}+\alpha\int_{\Si}u_\e^2dv_g+1+o_\d(1).
\ena
This together with (\ref{Gloc}), (\ref{Green-01}), and $c_\e u_\e \ra G_\alpha$ in $L^2(\Si,g)\cap C^1(\Sigma\setminus\{\mathbf{G}(x_0)\cup\mathbf{S}\})$ shows}
$$\label{ueintd}
\int_{\Si\setminus \bigcup_{k=1}^{I_0} \psi_k^{-1}(\mathbb{B}_{\delta})}|\nabla_g u_\e|^2 d v_g
=\f{1}{c_\e^2}\left(\f{1}{2\pi I_0}\log\d +A_0+\alpha\int_{\Si}G_\alpha^2dv_g+o_\e(1)+o_\d(1)\right).
$$
We then calculate
$$
\int_{ \bigcup_{k=1}^{I_0} \psi_k^{-1}(\mathbb{B}_{\delta})}|\nabla_g u_\e|^2 d v_g
=1-\f{1}{c_\e^2}\left(\f{1}{2\pi I_0}\log\d  +A_0+o_\e(1)+o_\d(1)\right):=\tau_\e.
$$
Set $s_\e=\sup_{\partial\psi^{-1}(\mathbb{B}_\d)}u_\e$ and $\hat{u}_\e=(u_\e-s_\e)^+.$ Clearly, $\hat{u}_\e\in W_0^{1,2}(\psi_k^{-1}(\mathbb{B}_\d))$. Moreover, we have
 $$
\int_{ \mathbb{B}_{\delta}}|\nabla_g (\hat{u}_\e\circ\phi^{-1})|^2 d x=\int_{ \psi_k^{-1}(\mathbb{B}_{\delta})}|\nabla_g \hat{u}_\e|^2 d v_g\leq\f{1}{I_0}\int_{ \bigcup_{k=1}^{I_0} \psi_k^{-1}(\mathbb{B}_{\delta})}|\nabla_g u_\e|^2 d v_g
\leq\f{\tau_\e}{I_0}.
$$
Then by using Lemma \ref{circle}, we obtain
\bea
{\nonumber}\limsup_{\e\ra0}\int_{ \psi^{-1}(\mathbb{B}_{\delta})}(e^{4\pi \ell\hat{u}_\e^2/\tau_\epsilon}-1) d v_g&=&
\limsup_{\e\ra0}\int_{\mathbb{B}_{\delta}}V(y)e^{2f}|y|^{2\beta}(e^{4\pi (1+\beta_0)I_0(\hat{u}_\e\circ\phi^{-1})^2/\tau_\epsilon} -1)d y \\[1.2ex]
{\nonumber}&=&\limsup_{\e\ra0}e^{o_\d(1)}\int_{\mathbb{B}_{\delta}}V_0|y|^{2\beta_0}(e^{4\pi (1+\beta_0)I_0(\hat{u}_\e\circ\phi^{-1})^2/
\tau_\epsilon} -1)dy\\[1.2ex]
\label{part1}&\leq& \f{\pi V_0e^{1+o_\d(1)}}{1+\beta_0}\d^{2+2\beta_0}.
\eea
For any fixed $R>0$, we have $u_\e/c_\e=1+o_\e(1)$  on $\psi_k^{-1}(\mathbb{B}_{Rr_\e^{1/(1+\beta)}})$ ($k=1,\cdots,I_0$). Hence, using the definition of $\tau_\e$, we obtain
\bna
(4\pi \ell-\e)u_\e^2&\leq& 4\pi \ell(\hat{u}_\e+s_\e)^2\\[1.2ex]
&=& 4\pi \ell \hat{u}_\e^2+8\pi \ell \hat{u}_\e s_\e+o_\e(1)\\[1.2ex]
&=& 4\pi\ell\hat{u}_\e^2-4(1+\beta_0)\log\d +8\pi \ell A_0+o(1)\\[1.2ex]
&=&  4\pi \ell_0 \hat{u}_\e^2/\tau_\epsilon-2(1+\beta_0)\log\d +4\pi \ell A_0+o(1),
\ena
where $o(1)\ra 0$ as $\epsilon\ra 0$ first, and then $\d\ra0$. Combining this with (\ref{part1}), we have
\bna
{\nonumber}\int_{\psi_k^{-1}(\mathbb{B}_{Rr_\e^{1/(1+\beta_0)}})}e^{(4\pi \ell-\e)u_\e^2} d v_g&\leq&
\d^{-2-2\beta} e^{4\pi \ell A_0+o(1)}\int_{\psi_k^{-1}(\mathbb{B}_{Rr_\e^{1/(1+\beta_0)}})}e^{(4\pi \ell-\e)\hat{u}_\e^2/
\tau_\epsilon} d v_g \\[1.2ex]
{\nonumber}&=&\d^{-2-2\beta_0} e^{4\pi \ell A_0+o(1)}\int_{\psi_k^{-1}(\mathbb{B}_{Rr_\e^{1/(1+\beta_0)}})}(e^{(4\pi \ell-\e)\hat{u}_\e^2
/\tau_\epsilon}-1) d v_g+o(1)\\[1.2ex]
{\nonumber}&\leq&\d^{-2-2\beta_0} e^{4\pi \ell A_0+o(1)}\int_{\psi_k^{-1}(\mathbb{B}_{\d})}(e^{(4\pi \ell-\e)\hat{u}_\e^2/
\tau_\epsilon}-1) d v_g+o(1)\\[1.2ex]
&\leq& \f{\pi V_0e^{1+4\pi \ell A_0+o(1)}}{1+\beta_0}.
\ena
Letting $\e\ra 0 $ first, and then $\d\ra 0$ , we obtain
\be\label{up1}
\limsup_{\e\ra0}\int_{\cup_{k=1}^{I_0}\psi_k^{-1}(\mathbb{B}_{Rr_\e^{1/(1+\beta_0)}})}e^{(4\pi \ell-\e)u_\e^2} d v_g\leq\f{\pi I_0V_0 e^{1+4\pi \ell A_0}}{1+\beta_0}.
\ee
Also we have
\bna\label{up2}
{\nonumber}\int_{\cup_{k=1}^{I_0}\psi_k^{-1}(\mathbb{B}_{Rr_\e^{1/(1+\beta_0)}})}e^{(4\pi \ell-\e )u_\e^2} d v_g&=&I_0(1+o_\e(1))\int_{\mathbb{B}_{R r_\e^{1/(1+\beta)}}}V_0e^{2f}|x|^{2\beta}\tilde{u}_\e^2 e^{(4\pi \ell-\e)\tilde{u}_{\e}^2}dx\\[1.2ex]
\label{up3}&=&\f{I_0\lambda_\e}{c_\e^2}(1+o_\e(1))\left(\int_{\mathbb{B}_R(0)}V_0|y|^{2\beta}e^{8\pi \ell \varphi_\e}dy+o_\e(1)\right)\\[1.2ex]
&=&\f{\lambda_\e}{c_\e^2}(1+o(1)),
\ena
where $o(1)\ra 0$ as $\epsilon\ra 0$ first, and then $R\ra\infty$. This together with (\ref{up1}) and (\ref{gama2}) leads to
 \bea
 {\nonumber}\sup_{u\in \mathscr{H}_{\mathbf{G}},\,\|u\|_{1,\alpha}\leq 1}\int_\Si e^{4\pi \ell u^2}dv_g &\leq&\mathrm{Vol}_g(\Si)+
 \lim_{R\ra \infty}\limsup_{\e\ra0}\int_{\cup_{k=1}^{I_0}\psi_k^{-1}(\mathbb{B}_{Rr_\e^{1/(1+\beta_0)}})}e^{(4\pi \ell-\e)u_\e^2} d v_g\\[1.2ex]
\label{upperb}&\leq&\mathrm{Vol}_g(\Si)+\f{\pi I_0V_0 e^{1+4\pi \ell A_0}}{1+\beta_0}.
 \eea

 \subsection{Existence of extremal functions.}

     Recall that $(\Si,g)$ has a conical singularity of the order $\beta_0$ at $x_0$  with $-1<\beta_0\leq 0$, $ I_0=I(x_0)$ and $ \beta_0=\beta(x_0)$, where $I(x)$ and $\beta(x)$ are defined as in (\ref{number}) and (\ref{betax}). In this section, we shall construct a sequence of functions $\widetilde{\Phi}_{\epsilon}\in \mathscr{H}_{\mathbf{G}}$ satisfying $\|\widetilde{\Phi}_{\e}\|_{1,\alpha}= 1$, and
\be\label{testbound}
 \int_{\Si} e^{4\pi\ell\widetilde{\Phi}_{\e}^2}dv_g >\mathrm{Vol}_g(\Si)+\f{\pi I_0V_0 e^{1+4\pi \ell A_0}}{1+\beta_0},
 \ee
  where $A_0$ and $V_0$ are constants defined as in (\ref{a-0}) and (\ref{V0}). The contradiction between (\ref{upperb}) and (\ref{testbound})
  implies that $c_\epsilon$ must be bounded, i.e. blow-up does not occur. This ends the proof of Theorem \ref{TH1}.

Set $R=(-\log\e)^{1/(1+\beta_0)}$.
 It follows that $R\ra\infty$ and $R\e\ra 0$ as $\e\ra0$. Hence, when $\e>0$ is sufficiently small, ${B}_{g_0,2R\e}(\sigma_i(x_0))\cup{B}_{g_0,2R\e}(\sigma_j(x_0))=\varnothing\,\textrm{for }1\leq i< j\leq I_0$. We firstly define a cut-off function $\eta$ on ${B}_{g_0,2R\e}(x_0)$, which is radially symmetric with respect to $x_0$. Besides, we require $\eta\in C^\infty_0( {B}_{g_0,2R\e}(x_0))$ to be a nonnegative function satisfying $\eta=1$ on ${B}_{g_0,R\e}(x_0)$ and $\|\nabla \eta\|_{L^\infty({B}_{2R\e})}
     =O(\f{1}{R\epsilon})$.
 Then we define a sequence of functions $\Phi_\e$ on $\Si$ for small $\e>0$ by
 \begin{small}
 \be\label{testPhi}\Phi_\e=\le\{
     \begin{array}{llll}
     &c+\dfrac{-\dfrac{1}{4\pi \ell}\log\le(1+\f{\pi I_0}{1+\beta_0}\frac{d_{g_0}(x,\sigma_k(x_0))^{2(1+\beta_0)}}{\epsilon^{2(1+\beta_0)}}\ri)+b}{c},
    \,&x\in \overline{B_{g_0,R\e}(\sigma_k(x_0))}\\[1.5ex]\\
    &\dfrac{G_\alpha(x)-\eta(\sigma_k^{-1}(x)) \widetilde{\Psi}(\sigma_k^{-1}(x))}{c},\quad &x\in {B}_{g_0,2R\e}(\sigma_k(x_0))\setminus \overline{{B}_{g_0,R\e}(\sigma_k(x_0))}\\[1.5ex]
     &\dfrac{G_\alpha}{c},\quad &x\in \Si\setminus \bigcup_{k=1}^{I_0}\sigma_k( {B}_{2R\e}),
     \end{array}
     \ri.
     \ee
 \end{small}
where $k$ is taken from $1$ to $I_0$, and $\widetilde{\Psi}$ is the function mentioned in (\ref{Gloc0}),
both $b$ and $c$ are constants depending on $\e$ to be determined later.

Recall (\ref{Green-01}), $G_\alpha(\sigma(x))=G_\alpha(x)$ for any $x\in {\Si}\setminus\bigcup_{k=1}^{I_0}\{\sigma_k(x_0)\}$ and all $\sigma\in\mathbf{G}$. Combining this with our premises that $\eta$ is radially symmetric and any $\sigma\in\mathbf{G}$ is a isometric map, for sufficiently small $\e$, we conclude
\be\label{phisym}
\Phi_\e(x)=\Phi_\e(\sigma(x))\quad \forall \sigma\in\mathbf{G},\quad \mathrm{a.e.}\,x\in\Si.
\ee
Set $\bar{\Phi}_\e=\Phi_\e-\f{1}{\mathrm{Vol}_g(\Si)}\int_{\Si}\Phi_\e dv_g$.
We shall choose suitable $b$ and $c$
to make $\widetilde{\Phi}_\e=\bar{\Phi}_\e/\|\bar{\Phi}_\e\|_{1,\alpha}\in \mathscr{H}_{\mathbf{G}}$.
Since the calculation is very similar to \cite[pages\,3365-3368]{2017JFA}, we omit the details but give its outline here.
{Integration by parts shows}
\begin{small}
\bna
\int_{\Si\setminus \bigcup_{k=1}^{I_0} B_{g_0,R\epsilon}(\sigma_k(x_0))}|\nabla_g G_\alpha|^2 d v_g
&=&- \sum_{k=1}^{I_0}\int_{\partial B_{g_0,R\epsilon}(\sigma_k(x_0))}G_\alpha\frac{\p G_\alpha}{\p n} ds_g+\int_{\Si\setminus \bigcup_{k=1}^{I_0} B_{g_0,R\epsilon}(\sigma_k(x_0))}G_\alpha\Delta_g G_\alpha d v_g\\[1.2ex]
&\,&=-\f{1}{2\pi I_0}\log R\e  +A_0+\alpha\int_{\Si}G_\alpha^2dv_g+O(R\e),
\ena
\end{small}
and it follows that
$$
\int_{\Si\setminus \bigcup_{k=1}^{I_0} B_{g_0,R\epsilon}(\sigma_k(x_0))}|\nabla_g \Phi_\e|^2 d v_g=\f{1}{c^2}\le(-\f{1}{2\pi I_0}\log R\e  +A_0+\alpha\int_{\Si}G_\alpha^2dv_g+O(R\e)\ri).
$$
Here we use estimates
$$
\int_{B_{g_0,2R\epsilon}(\sigma_k(x_0))\setminus B_{g_0,R\epsilon}(\sigma_k(x_0))}|\nabla_g (\widetilde{\Psi}\eta)|^2 d v_g=O(R^2\e^2),
$$
and
$$\int_{B_{g_0,2R\epsilon}(\sigma_k(x_0))\setminus B_{g_0,R\epsilon}(\sigma_k(x_0))}\nabla_g G_\alpha\nabla_g (\widetilde{\Psi}\eta) d v_g=O(R\e).
$$
By a straightforward calculation, we obtain
$$
\int_{\bigcup_{k=1}^{I_0}B_{g_0,R\epsilon}(\sigma_k(x_0))}|\nabla_g \Phi_\e|^2 d v_g=\f{1}{4\pi\ell c^2}\le(\log\f{\pi I_0}{1+\beta_0}-1+\log R^{2+2\beta_0}+O(R^{-2-2\beta_0})\ri).
$$
Thus
\bna
\int_{\Si}|\nabla_g \Phi_\e|^2 d v_g&=&\f{1}{c^2}\le(-\f{\log\e}{2\pi I_0}+A_0+\alpha\int_{\Si}G^2dv_g-\f{1}{4\pi\ell}\ri.\\[1.2ex]
&\;&+\le.\f{1}{4\pi\ell}\log\f{\pi I_0}{1+\beta_0}+O(R^{-2-2\beta_0}) \ri).
\ena
Moreover, we have
$$
\int_{\Si}|\Phi_\e-\bar{\Phi}_\e|^2 dv_g=\f{1}{c^2}\le(\int_{\Si}G^2dv_g+O(R^{-2-2\beta_0})\ri).
$$
In the view of $\widetilde{\Phi}_\e\in W^{1,2}(\Si,g)$ and $\|\widetilde{\Phi}_\e\|_{1,\alpha}=1$, it follows from the above equations that
$$
c^2=-\f{1}{2\pi I_0}\log \e+A_{0}-\f{1}{4\pi \ell}+\f{1}{4\pi \ell}\log \f{\pi I_0}{1+\beta_0}+O(R^{-2(1+\beta_0)}),
$$
and
$$
b=\f{1}{4\pi \ell}+O(R^{-2(1+\beta_0)}).
$$
On $B_{g_0,R\epsilon}(\sigma_k(x_0))$, we have the following estimate:
     \bna
     4\pi\ell(1+\beta_0)\widetilde{\Phi}_\e^2&\geq& -2\log\le(1+\f{\pi I_0 }{1+\beta_0}\frac{r^{2(1+\beta_0)}}{\epsilon^{2(1+\beta_0)}}\ri)+1-2(1+\beta_0)\log\e\\[1.2ex]
     &\,&+4\pi\ell A_0+\log\f{\pi I_0}{1+\beta_0}+O(R^{-2-2\beta_0}).
     \ena
 Note that
 $$
 \int_{\mathbb{B}_R}\dfrac{1}{\le(1+\f{\pi I_0 }{1+\beta_0}{|y|^{2(1+\beta_0)}}\ri)|y|^{2\beta_0}}=1-\dfrac{1}{1+\f{\pi I_0}{1+\beta_0}R^{2+2\beta_0}}.
 $$
 This leads to
\bna
\int_{\bigcup_{k=1}^{I_0}B_{g_0,R\epsilon}(\sigma_k(x_0))}{e^{(4\pi \ell-\e)\widetilde{\Phi}_\e^2}} dv_g&=&\left(1+O(R\e)\right)I_0\int_{\mathbb{B}_{R\e}}V_0|x|^{2\beta_0}{e^{(4\pi \ell-\e)(\widetilde{\Phi}_\e^2(\exp_{x_0}x)^2}} dx\\[1.2ex]
&\geq&\left(1+O(R\e)\right)\f{\pi V_0 I_0e^{1+4\pi \ell A_0}}{1+\beta_0}.
\ena
On the other hand, by $e^{(4\pi \ell-\e)\widetilde{\Phi}_\e^2}\geq1+(4\pi\ell-\e)\widetilde{\Phi}_\e^2$, we obtain
$$
 \int_{\Si\setminus \bigcup_{k=1}^{I_0} B_{g_0,2R\epsilon}(\sigma_k(x_0)}{e^{(4\pi \ell-\e)\widetilde{\Phi}_\e^2}} dv_g
 \geq \mathrm{Vol}_g(\Si)+\f{4\pi \ell}{c^2}\int_{\Si}G^2dv_g+ O(R^{-2-2\beta_0}),
$$
which immediately lead to
\bna
\int_{\Si}{e^{(4\pi \ell-\e)\widetilde{\Phi}_\e^2}} dv_g&=& \int_{\cup_{k=1}^{I_0} B_{g_0,2R\epsilon}(\sigma_k(x_0))}e^{(4\pi \ell-\e)\widetilde{\Phi}_\e^2} dv_g+\int_{\Si\setminus\cup_{k=1}^{I_0} B_{g_0,2R\epsilon}(\sigma_k(x_0))}e^{(4\pi \ell-\e)\widetilde{\Phi}_\e^2} dv_g\\[1.2ex]
&\geq&\mathrm{Vol}_g(\Si)+\f{\pi V_0 I_0e^{1+4\pi \ell A_0}}{1+\beta_0}+\f{4\pi \ell}{c^2}\int_{\Si}G^2dv_g+ O(R^{-2-2\beta_0}).
\ena
Note that $R=(-\log\e)^{1/(1+\beta_0)}$, and $O(R^{-2(1+\beta_0)})=o(1/c^2)$. If $\epsilon>0$ is chosen sufficiently small,
then we arrive at (\ref{testbound}), as desired. $\hfill\Box$
\section{Proof of Theorem \ref{TH2}. }\label{pf-2}
The method we use to proof of Theorem \ref{TH2} is analogous to that of Theorem \ref{TH1}. Firstly, we conclude $4\pi \ell$ is the best constant for the inequality (\ref{Th2}) by a discussion totally similar to that in Subsection 2.1.
 Then we introduce an orthonormal basis  $(e_{j})$ $(1\leq j\leq n_\ell)$ of $E_{\ell}$ satisfying
$$
\left\{
\begin{array}{lllllllll}
E_{\ell}=\mathrm{span}\{e_{1}\;,\cdots,\;e_{n_{\ell}}\},\quad\\[1.5ex]
e_j\in C^{0}(\Sigma,g)\cap\mathscr{H}_{\mathbf{G}}, &\forall 1\leq j\leq n_\ell\\[1.5ex]
\int_{\Si}|e_{j}|^{2}dv_g=1,&\forall 1\leq j\leq n_\ell\\[1.5ex]
\int_{\Si}e_{l}e_{m}dv_g=0,& m\not=l,
\end{array}
\right.
$$
where $n_\ell=\mathrm{dim}\,E_{\ell}$. Under this orthonormal basis, $E_{\ell}^{\bot}$ is written as
$$\label{E}
 E_{\ell}^{\bot}=\bigg{\{}u\in \mathscr{H}_{\mathbf{G}}:\int_{\Si}ue_{j}dv_g=0,\;1\leq j\leq n_\ell\bigg{\}}.
$$
Secondly, we prove the existence of extremals for subcritical Trudinger-Moser functionals. Namely, for any $0<\epsilon<4\pi\ell$, there exists some $u_{\epsilon}\in E_{\ell}^{\bot}\cap  C^{1}(\Sigma\setminus\{p_1,\cdots,p_L\},g_0)\cap C^0(\Sigma,g_0)$  such that
$$\int_{\Si}e^{(4\pi\ell-\epsilon) u_{\epsilon}^{2}} {d}v_g=\sup\limits_{ u\in E_{\ell}^{\bot}\,\|u\|_{1,\alpha}\leq 1}\int_{\Si}e^{(4\pi\ell-\epsilon) u^{2}}dv_g.$$
Clearly $u_{\epsilon}$ satisfies the Euler-Lagrange equation
\begin{equation}\label{l-E=L}
\left\{
\begin{array}{lllllllll}
\triangle_g u_{\epsilon} -\alpha u_{\epsilon}=\frac{1}{\lambda_{\epsilon}}u_{\epsilon}e^{(4\pi\ell-\epsilon) u_{\epsilon}^{2}}-\f{\mu_\epsilon}{\lambda_\epsilon}-\sum_{j=1}^{n_\ell}\omega_{j,\e}e_j,
\\[1.5ex]
\|u_\e\|_{1,\alpha}= 1,\\[1.5ex]
\lambda_\epsilon=\int_\Si u_\epsilon^2 e^{(4\pi \ell-\epsilon)u_\epsilon^2}dv_g,\\[1.5ex]
\mu_\epsilon=\f{1}{\mathrm{Vol}_g(\Si)}\int_\Si u_\epsilon e^{(4\pi \ell-\epsilon)u_\epsilon^2}dv_g,\\[1.5ex]
\omega_{j,\e}=\f{1}{\lambda_\epsilon}\int_\Si e_ju_\epsilon e^{(4\pi \ell-\epsilon)u_\epsilon^2}dv_g.
\end{array}
\right.
\end{equation}
Assume $u_\e$ converges to $u^*$ weakly in $W^{1,2}(\Sigma,g)$,
   strongly in $L^s(\Sigma,g)$ for any $s>1$, and almost everywhere in $\Sigma$.
    If $u_{\epsilon}$ is uniformly bounded, then we have by the Lebesgue dominated convergence theorem
$$\int_{\Si}u^*e_{j} dv_g=\lim\limits_{\epsilon\rightarrow0}\int_{\Si}u_\e e_{j} dv_g=0,\quad\forall 1\leq j\leq n_\ell,$$
and thus $u^*\in  E_{\ell}^{\bot}\cap  C^{1}(\Sigma\setminus\{p_1,\cdots,p_L\},g_0)\cap C^0(\Sigma,g_0)$ is the desired extremal function.

If blow-up happens, we still have analogs of Lemmas \ref{blow-upanal} and \ref{beta=0}. For any $1<q<2$, we obtain $c_\e u_\e\rightharpoonup G$ weakly in $W^{1,q}(\Si,g)$ , where $G$ is green function satisfying
$$
\le\{
  \begin{array}{lll}
\Delta_g G -\alpha G=\sum_{i=1}^{I_0}\f{\delta_{\sigma_i(x_0)}}{I_0}-\f{1}{\mathrm{Vol_g(\Si)}}-\sum_{j=1}^{n_\ell}e_j(x_0)e_j\\[1.2ex]
\int_{\Si} Ge_j dv_g=0,\quad  1\leq j\leq n_\ell\\[1.2ex]
G(\sigma_i(x))=G(x),\quad 1\leq i\leq N,\,x\in\Sigma\setminus\{\sigma_i(x_0)\}_{i=1}^{I_0}.
\end{array}
  \ri.
$$
Similar to the proof of (\ref{upperb}), we can drew a conclusion that
\be\label{upperbl}
\sup_{u\in E_{\ell}^{\bot},\,\|u\|_{1,\alpha}\leq 1}\int_{\Si}e^{4\pi \ell u^2}dv_g\leq \mathrm{Vol}_g(\Si)+\f{\pi I_0V_0}{1+\beta_0}e^{1+4\pi \ell A_0},
\ee
where all the constants in (\ref{upperbl}) have the same definition as in the last section.

At last, we shall construct a sequence of functions to contradict (\ref{upperbl}). Denote
$$
{\omega_{\epsilon}}={\Phi}_{\epsilon}-\sum\limits_{j=1}^{n_{j}}({\Phi}_\e,e_{j})e_j,
$$
where $\Phi_\e$ is defined as in (\ref{testPhi}), and
$$
({\Phi}_\e,e_{j})=\int_{\Si}{\Phi}_{\epsilon}e_{j} dv_g .
$$
Set $\widetilde{\omega}_\e=\omega_\e-\f{1}{\mathrm{Vol}_g(\Si)}\int_{\Si}\omega_\e dv_g$. We may choose suitable constants $b$ and $c$ to make $\widetilde{\omega}_\e\in E_{\ell}^{\bot}$. A straightforward calculation shows
\begin{equation*}
\begin{split}
\int_{\Si}e^{4\pi\ell\f{{\widetilde{\omega}_{\epsilon}}}{\|{\widetilde{\omega}_{\epsilon}}\|_{1,\alpha}}}dv_g=&\int_{\Si}e^{4\pi\ell {\widetilde{\omega}_{\epsilon}}^{2}+o\big{(}\frac{1}{\log{\epsilon}}\big{)}}dv_g\\
&\geq \big{(}1+o\big{(}\dfrac{1}{\log{\epsilon}}\big{)} \big{)} \big{(}\mathrm{Vol}_g(\Si)+4\pi I_0 \dfrac{\|G\|^{2}_{2}}{c^{2}}+\f{\pi I_0V_0 e^{1+4\pi \ell A_0}}{1+\beta_0}\big{)} \\
&\geq\mathrm{Vol}_g(\Si)+4\pi I_0 \dfrac{\|G\|^{2}_{2}}{-\log\e}+\f{\pi I_0V_0 e^{1+4\pi \ell A_0}}{1+\beta_0}+o\big{(}\frac{1}{\log{\epsilon}}\big{)}.\\
\end{split}
\end{equation*}
This  indicates
$$
\sup_{u\in E_{\ell}^{\bot},\,\|u\|_{1,\alpha}\leq 1}\int_{\Si}e^{4\pi \ell u^2}dv_g>\mathrm{Vol}_g(\Si)+\f{\pi I_0V_0}{1+\beta_0}e^{1+4\pi \ell A_0},
$$
which  contradicts (\ref{upperbl}). Thus the proof of Theorem \ref{TH2} is finished.$\hfill\Box$

\end{document}